\documentclass[aip,reprint,jcp,nofootinbib]{revtex4-1}
\usepackage{graphicx}
\usepackage{dcolumn}
\usepackage{bm}
\usepackage{amssymb}
\usepackage{amsmath}
\usepackage[usenames]{color} 
\usepackage{ulem} 
\usepackage{epstopdf}
\usepackage{latexsym}
\usepackage{float}

\newcommand{\half}{{\tfrac{1}{2}}}

\begin{document}

\preprint{AIP/123-QED}

\title[]{Reliability Assessment for Large-Scale
  Molecular Dynamics Approximations}

\author{Francesca Grogan}
\email{fgrogan@ucsd.edu}
\affiliation{Center for Computational Mathematics, Department
  of Mathematics, University of California at San Diego, La Jolla, CA 92093, USA}
\author{Michael  Holst}
\email{mholst@ucsd.edu}
\author{Lee Lindblom }
\email{llindblom@ucsd.edu}
\affiliation{Center for Computational Mathematics, Department
  of Mathematics, University of California at San Diego, La Jolla, CA 92093, USA}
\affiliation{ Center for
  Astrophysics and Space Sciences, University of California at San
  Diego, La Jolla, CA 92093, USA}
\author{Rommie Amaro}
\email{ramaro@ucsd.edu}
\affiliation{Department of Chemistry and Biochemistry,
  University of California at San
  Diego, La Jolla, CA 92093, USA}
  
\date{\today}

\begin{abstract}
Molecular dynamics (MD) simulations are used in biochemistry, physics,
and other fields to study the motions, thermodynamic properties, and
the interactions between molecules.  Computational limitations and the
complexity of these problems, however, create the need for
approximations to the standard MD methods and for uncertainty
quantification and reliability assessment of those approximations.  In
this paper, we exploit the intrinsic two-scale nature of MD to
construct a class of large-scale dynamics approximations.  The
reliability of these methods are evaluated here by measuring the
differences between full, classical MD simulations and those based on
these large-scale approximations.  Molecular dynamics evolutions are
non-linear and chaotic, so the complete details of molecular
evolutions cannot be accurately predicted even using full, classical
MD simulations.  This paper provides numerical results that
demonstrate the existence of computationally efficient large-scale MD
approximations which accurately model certain large-scale properties
of the molecules: the energy, the linear and angular momenta, and
other macroscopic features of molecular motions.
\end{abstract}

\pacs{83.10.Mj, 83.10.Rs, 87.10.Tf, 87.15.ap}

\keywords{molecular dynamics, uncertainty quantification, numerical
  methods, normal modes}
\maketitle

\section{Introduction}
\label{s:Intro}

Traditional molecular dynamics (MD) simulations use Newton's classical
equations of motion, with an effective potential that models the
interactions between atoms, to describe the evolution of
molecules.~\cite{McCammon1987, Allen1987, Frenkel2002} This standard
MD method has been applied to a variety of problems in biochemistry
and condensed matter physics in recent years.~\cite{Karplus2002,
  LASKOWSKI20091692, Durrant2011, Li2011aom, Wang2011,
  Steinhauser2012, Perilla2015, Vettoretti2016}

The results obtained from these MD simulations can be characterized,
qualitatively, as the evolution of the overall position and
orientation of each molecule, plus vibrations of the individual atoms
about their average positions within each molecule.  The timescales
associated with the macroscopic bulk motions of the molecules are
typically much longer than the timescales associated with the internal
atomic vibrations.  However, these very short timescale vibrations
determine the maximum timesteps allowed for accurate solutions of the
classical MD evolution equations using standard numerical methods.
This fundamental fact, together with the need to simulate extremely
large complex molecules in modern biochemical research, means the
computational cost of performing MD simulations can be prohibitively
large.  Computational cost is therefore one of the factors that drives
the need to develop approximation methods capable of obtaining
reliable simulations of those aspects of the molecular systems of
primary interest to researchers.

Another factor that motivates the development of approximation methods
for MD simulations is the well-known fact that typical $n$-body
systems exhibit chaotic behavior in which exponentially divergent
evolutions result from small perturbations of initial
conditions.\cite{Poincare1890} It is simply impossible, and therefore
pointless to attempt, to simulate all the details in the evolutions of
complex molecular systems.  Many large-scale properties of such
systems, including statistical time averages of various properties,
are nevertheless well defined~\cite{Allen1987, Bond2007,Leim2004} and
these properties are therefore in principle observable and
simulatable.  These macroscopic characteristics include the positions,
average velocities, and other time averages of thermodynamic
quantities.

Given these two motivating factors---the need for greater
computational efficiency and the fundamental inability to simulate
chaotic dynamics in complete microscopic detail---we have developed a
new class of large scale MD approximations.  We construct these
approximations by starting with a new representation of the exact
traditional MD equations of motion.  The standard MD equations
generally use the Cartesian coordinates of the location of each atom
as the primary dynamical variables.  We transform these Cartesian
variables into a new representation that separates the macroscopic
location and orientation of each molecule, from the internal degrees
of freedom that represent the molecular vibrations.  This
transformation is (loosely) motivated by Wilson's representation of
the internal degrees of freedom of a molecule by normal
modes.\cite{Wilson1941} Our approach differs from Wilson's, however,
by providing an exact representation of the molecular motions even
when the mode amplitudes are not small.  We test this new mode-basis
description of MD by comparing the results of numerical simulations
using it with those based on the standard classical MD equations.  The
results of these numerical tests, described in more detail in
Sec.~\ref{s:TestingModeBasisRepresentation}, confirm that our new
mode-basis representation of MD is exact.

In our mode-basis representation of MD, the large-scale degrees of
freedom are cleanly separated from the collection of mode amplitudes
that describe the internal vibration degrees of freedom of the
molecule.  It is straightforward to replace the exact equations that
determine the evolution of these mode amplitudes, with various
approximate expressions.  For example the mode amplitudes could be
chosen to satisfy the analytic sinusoidal-in-time expressions derived
from the small-amplitude normal-mode equations.  Another possibility
would be to use the normal-mode sinusoidal-in-time expressions for
some number of the highest frequency modes, while evolving the
remaining lower frequency mode amplitudes numerically using the exact
evolution equations.  This approach would be qualitatively similar to
that used in the constraint algorithms like SHAKE and
RATTLE.\cite{Ryckaert1977,Andersen1983} Another possibility would
simply be to ignore the internal degrees of freedom of the atoms
completely by setting the mode amplitudes equal to zero, i.e. to their
expected long timescale averages.  We present numerical tests that
compare different approximations of this type with the results of the
exact MD evolutions.  The results of these numerical tests are
described in more detail in
Sec.~\ref{s:TestingLargeScaleApproximations}.

We assess the reliability of our new large-scale approximation methods
by applying the techniques of uncertainty quantification (UQ).\cite{
  Smith2014, Patrone2016, Reeve2017, Rizzi2012, Rizzi2_2012} Here we
focus attention on measuring the accuracy of our new large-scale
approximation methods by comparing them to standard classical MD,
setting aside other important issues such as time-stepper integration
errors, errors in the molecular interaction potential model, errors
that arise from the use of classical rather than a fully quantum
description of MD, etc.  The MD evolution equations are highly
non-linear, as are the equations for our large-scale approximations.
It is possible to derive rigorous analytic mathematical bounds on the
errors in our large-scale approximation methods.\cite{Grogan2017}
However the bounds we have obtained in this way are quite weak, and do
not provide a good estimate of the size of the actual large-scale
approximation errors in practical simulations.  Therefore we focus the
discussion of our UQ analysis in this paper on making detailed
numerical comparisons between full classical MD simulations and those
obtained for identical molecular systems using our new large-scale
approximations.

Since MD simulations are typically performed to estimate the values of
various macroscopic observables of the molecules, we have focused our
uncertainty quantification analysis on assessing the errors in the
large-scale approximation values of those quantities.  In particular
we evaluate the errors in the energy, the linear and angular momenta,
and the errors in the positions and orientations of each molecule.
Our numerical results show that all the large-scale approximations
tested here are linear momentum conserving, and consequently the
positions and velocities of the molecular centers of mass are also
determined exactly.  Some of our large-scale approximations tested
here also conserve energy and angular momentum exactly.  Angular
momentum conservation does not, however, guarantee that the
orientations or angular velocities of the molecules are determined
accurately.  We show that these orientation features evolve
chaotically in MD systems, and are therefore unpredictable even in
full classical MD simulations.

The remainder of this paper is organized as follows.  In
Sec.~\ref{s:ClassicalMD}, we derive a new mode-basis representation of
the molecular dynamics evolution equations, and then use them to
derive several large-scale molecular dynamics approximations.  We have
implemented these equations in a numerical MD evolution code and have
used this code to evolve simple models of several of the smaller
fullerene molecules: $C_{20}$, $C_{26}$, $C_{60}$ and $C_{70}$.  In
Sec.~\ref{s:ReliabilityAssessment} we discuss the results of our
numerical simulations of these molecules using both standard classical
MD and several large-scale approximations, highlighting the
uncertainty quantification of the macroscopic properties of these
molecules.  We conclude by summarizing and discussing our results
briefly in Sec.~\ref{s:Discussion}.

\section{Classical Molecular Dynamics}
\label{s:ClassicalMD}

Classical molecular dynamics (MD) uses Newton's equations to describe
the motions of the atoms (represented as point particles) that make up
the molecules being studied.  We use the notation $\vec x_A(t)$ to
denote the Cartesian coordinates of the location of atom $\scriptstyle
A$ as a function of time.  The classical MD equations of motion for
the atoms that make up a molecule (or collection of molecules) are
therefore given by

\begin{eqnarray}
m_A \frac{d^{\,2} \vec x_A}{dt^2}&=&
-\frac{\partial U(\vec x_B)}{\partial \vec x_A},
\label{e:NewtonsLaw}
\end{eqnarray}
where $m_A$ is the mass of atom $\scriptstyle A$ and $U(\vec x_B)$ is
the effective potential energy function that describes the
interactions between the atoms.  This potential energy will in general
be a non-linear function of the locations of all of the atoms.
Consequently the MD evolution equations, Eq.~(\ref{e:NewtonsLaw}), are
non-linear and strongly coupled.  The solutions to these equations
therefore display the typical characteristics of chaotic $n$-body
systems,\cite{Allen1987, Bond2007,Leim2004} i.e. while most details
of the evolution of a particular initial molecular state cannot be
predicted, certain ``statistical'' features of the evolution can.

In this section we develop a new class of large-scale approximations
to the classical MD equations of motion, Eq.~(\ref{e:NewtonsLaw}).
These new approximations are designed to provide a more efficient way
to evaluate some of the observable ``statistical'' or
``thermodynamic'' features of MD systems.  We construct these
approximations in two steps.  In the first step, in
Sec.~\ref{s:ModeBasisMD}, we transform the exact classical MD
equations into a representation that cleanly separates some of the
observable macroscopic degrees of freedom of a molecule from the
internal degrees of freedom that determine its microscopic state. We
refer to this new representation of the MD equations as a mode-basis
representation, because the choice of variables used to describe the
microscopic state of a molecule is based (loosely) on the normal-mode
description of molecular vibrations.\cite{Wilson1941,Levitt1985} The
transformation we use to construct this mode-basis representation is
exact, however, so it is simply a change of variables for the
classical MD system given in Eq.~(\ref{e:NewtonsLaw}).  We present the
exact MD equations of motion for the mode-basis representation in
Sec.~\ref{s:md_large-scale}.  The second step in the development of
the new class of MD approximations is to replace the exact
mode-amplitude evolution equations with various approximations.  We
present several examples in Sec.~\ref{s:LargeScaleMD} that use simple
analytic expressions to approximate the evolutions of some or all of
the mode amplitudes which describe the internal molecular degrees of
freedom.

\subsection{Mode-basis Representation of MD}
\label{s:ModeBasisMD}

The mode-basis representation of MD is obtained by transforming the
Cartesian coordinate variables, $\vec x_A$, used in the standard
representation into new variables that separate into {\it i)} a set
that describes the macroscopic location and orientation of the
molecule, and {\it ii)} another set that describes the molecules'
microscopic vibrational dynamics.  The macroscopic location of a
molecule is represented by its center of mass, $\vec x_{CM}(t)$,
defined by
\begin{eqnarray}
\vec x_{CM}(t) = \frac{1}{M}\sum_A m_A \vec x_A(t),
\label{e:CMDef}
\end{eqnarray}
where $M=\sum_A m_A$ is the total mass of the molecule.
\footnote{We note that for
  MD simulations of collections of molecules, the single index
  $\scriptstyle A$ that identifies the atoms should be replaced by a
  pair of indices $\scriptstyle{mA}$, with $\scriptstyle m$ identifying
  the particular molecule and $\scriptstyle A$ the atom within that
  molecule.  All the macroscopic properties of these molecules should
  also acquire an additional $\scriptstyle m$ index to identify which
  molecule the attribute belongs to, e.g. $M$ should become $M_m$,
  $\vec x_{CM}(t)$ should become $\vec x_{CMm}(t)$, etc.  For
  simplicity of notation we will suppress these molecule-identifying
  indices in the discussion in this paper.}
The macroscopic
orientation of a molecule is represented by a time-dependent rotation
matrix ${\mathbf R}(t)$.  This matrix provides the transformation
between a reference frame fixed to and co-moving with the molecule,
and the global inertial frame used to describe the atoms in the
standard representation of MD.  We use the notation $\vec x_{0A} +
\delta \vec x_A(t)$ to denote the location of atom $\scriptstyle A$ in
the molecule's co-moving reference frame, where $\vec x_{0A}$
represents the time independent equilibrium position, and $\delta \vec
x_A(t)$ the displacement from equilibrium (which is not assumed to be
small) of this atom.  The global Cartesian coordinate location of atom
$\scriptstyle A$ is determined by these macroscopic variables---the
center of mass, $\vec x_{CM}(t)$ and the orientation matrix, ${\mathbf
  R}(t)$---along with the internal dynamical variables $\delta \vec
x_A$:
\begin{eqnarray}
  \vec x_A(t) = \vec x_{CM}(t) + {\mathbf R}(t)\cdot
  \bigl[\vec x_{0A} + \delta \vec x_A(t)\bigr].
  \label{e:ModeBasisTransformation0}
\end{eqnarray}

In the mode-basis representation of MD, the internal microscopic
degrees of freedom of a molecule are described by the variables
$\delta \vec x_A$.  Unfortunately these variables are not independent,
so special care must be taken to isolate the truly independent
degrees of freedom they represent.  To see this more clearly, note
that in the standard MD representation there are $3N$ variables, $\vec
x_A$, needed to represent the configuration state of a molecule having
$N$ atoms.  The macroscopic variables introduced above, $\{\vec
x_{CM}, {\mathbf R}\}$, represent $6$ of these degrees of freedom
(since there is a three-dimensional space of rotation matrices
$\mathbf R$).  Consequently there can only be $3N-6$ truly
independent internal microscopic degrees of freedom among the $3N$
variables $\delta \vec x_A$.  To isolate these independent degrees of
freedom we introduce a collection of ``mode-basis'' vectors $\vec
e^{\,\mu}_A$, where the index $\mu$ labels the $3N-6$ vectors
representing those independent degrees of freedom.  Without loss of
generality we can normalize these basis vectors:
\begin{equation}
\delta^{\mu\nu} = \sum_A \frac{m_A}{M}\, \vec e^{\,\mu}_A \cdot\vec
e^{\,\nu}_A,
\label{e:ModeOrthogonality}
\end{equation}
where $\delta^{\mu\nu}$ is the Kronecker delta.  To ensure that the
$\vec e^{\,\mu}_A$ are independent from the macroscopic variables,
$\{\vec x_{CM}, {\mathbf R}\}$, we choose them to be orthogonal to any
overall translation or rotation of the molecule:
\begin{eqnarray}
  0&=&\sum_A \frac{m_A}{M} \vec e^{\,\mu}_A,\label{e:TranslationConstraint}\\ 
  0&=&\sum_A \frac{m_A}{M} \vec e^{\,\mu}_A\times \vec x_{0A}.
  \label{e:RotationConstraint}
\end{eqnarray}
Appendix~\ref{s:NormalModeBasis} explains in detail why
Eqs.~(\ref{e:TranslationConstraint}) and (\ref{e:RotationConstraint})
are the conditions needed to enforce the translation and rotation
invariance of the eigenvectors $\vec e^{\,\mu}_A$.  Given any
collection of mode-basis vectors satisfying
Eqs.~(\ref{e:ModeOrthogonality})--( \ref{e:RotationConstraint}), it is
straightforward to write down a general expression for $\delta \vec
x_A$ in terms of $3N-6$ independent mode-amplitude functions
${\mathcal A}_\mu(t)$:
\begin{eqnarray}
  \delta \vec x_A&=& \sum_\mu {\mathcal A}_\mu(t)\, \vec e^{\,\mu}_A.
  \label{e:ModeSumd}
\end{eqnarray}
Using this expression and Eq.~(\ref{e:ModeBasisTransformation0}), it
is now possible to write down the transformation between the
mode-basis representation variables, $\{\vec x_{CM},{\mathbf
  R},{\mathcal A}_\mu\}$, and the Cartesian coordinates, $\vec x_A$
used in the standard representation of MD:
\begin{eqnarray}
  \vec x_A(t) = \vec x_{CM}(t) + {\mathbf R}(t)\cdot
  \left[\vec x_{0A} + \sum_\mu{\mathcal A}_\mu(t)\,\vec e^{\,\mu}_A\right].
  \label{e:ModeBasisTransformation1}
\end{eqnarray}

There are an infinite number of ways to choose the mode-basis vectors,
$\vec e^{\,\mu}_A$.  One natural choice is to let the $\vec
e^{\,\mu}_A$ be the eigenvectors of the Hessian of the potential
energy function:
\begin{equation}
0=-m_A\omega_\mu^2\,\vec e^{\,\mu}_A + \sum_B 
\frac{\partial^2 U}{\partial {\vec x_B}\partial {\vec x_A}}\cdot \vec e^{\,\mu}_B.
\label{e:ModeEquation}
\end{equation}
The Hessian matrix, $\partial^2 U/\partial {\vec x_B}\partial {\vec
  x_A}$, in this equation is to be evaluated at the equilibrium state
of the molecule where $\vec x_A = \vec x_{0A}$.  Since the Hessian is
a symmetric real matrix, the eigenvalues $\omega_\mu^2$ and
eigenvectors $\vec e^{\,\mu}_A$ are also real.
Equation~(\ref{e:ModeEquation}) is equivalent to Newton's equation of
motion, Eq.~(\ref{e:NewtonsLaw}), for the case of very small amplitude
oscillations about its equilibrium state.  This choice of eigenvectors
is therefore particularly useful for isolating the individual
microscopic degrees of freedom of a molecule.  The connection of these
eigenvectors to the classical normal-mode analysis of molecular
vibrations\cite{Wilson1941,Levitt1985} motivated our choice of the
name ``mode-basis'' representation.  The eigenvalues $\omega_\mu^2$ in
Eq.~(\ref{e:ModeEquation}) are non-negative for any stable molecule,
and except for the six zero-frequency modes that correspond to rigid
rotations and translations of the molecule, these eigenvalues are
strictly positive (generically).  Since the Hessian matrix is
symmetric, the eigenvectors $\vec e^{\,\mu}_A$ form a complete basis
for the $\delta \vec x_A$ that satisfy (or in the case of degenerate
eigenvalues, can be chosen to satisfy) the orthogonality conditions,
Eq.~(\ref{e:ModeOrthogonality}).  Appendix~\ref{s:NormalModeBasis}
shows that they also satisfy the constraints,
Eqs.~(\ref{e:TranslationConstraint}) and (\ref{e:RotationConstraint}).

The standard MD equations of motion, Eq.~(\ref{e:NewtonsLaw}), are
second-order ordinary differential equations.  Therefore both the
position $\vec x_A$ and the velocity $\vec v_A=d\vec x_A/dt$ of each
atom are needed to determine the full dynamical state of a molecule.
The analog of these velocity degrees of freedom for the macroscopic
variables are the center of mass velocity $\vec v_{CM}=d\vec
x_{CM}/dt$, and the time derivative of the orientation matrix
$d{\mathbf R}/dt$.  It is convenient and customary to express the time
derivative of the orientation matrix as an angular velocity vector,
$\vec\Omega$, in the following way.  The matrix ${\mathbf{*\Omega}}$
defined by,
\begin{equation}
\mathbf{*\Omega}=-\frac{d\mathbf{R}}{dt}\cdot \mathbf{R}^{-1},
\label{e:AngularVelocityMatrix}
\end{equation}
is anti-symmetric for any rotation matrix $\mathbf R$.  Therefore
${\mathbf{*\Omega}}$ is dual to a vector $\vec \Omega$:
\begin{equation}
*\Omega_{ij} = \sum_k\epsilon_{ijk}\Omega^k,
\label{e:AngularVelocityDef}
\end{equation}
where $\epsilon_{ijk}$ is the totally anti-symmetric tensor with
$\epsilon_{xyz}=1$ in Cartesian coordinates.  The time derivative
$d{\mathbf R}/dt$ is therefore given by,
\begin{eqnarray}
\frac{d\mathbf{R}}{dt}=-\mathbf{*\Omega}\cdot\mathbf{R},
\label{e:dRdt}
\end{eqnarray}
or in component notation,
\begin{eqnarray}
\frac{dR_{ij}}{dt}=-\sum_{k\ell}\epsilon_{ik\ell}\,\Omega^\ell R_{kj}.
\label{e:dRdtcomponents}
\end{eqnarray}

Using these results we can now write down the complete set of
transformation equations between the Cartesian coordinate variables,
$\{\vec x_A,\vec v_A\}$, used in the standard representations of
classical MD, and the variables, $\{\vec x_{CM},\vec v_{CM},{\mathbf
  R}, \vec \Omega, {\mathcal A}_\mu,d{\mathcal A}_\mu/dt\}$ used in
our new mode-basis representation:
\begin{eqnarray}
\vec x_A(t)&=& \vec x_{CM}(t) + \Delta \vec x_A(t),
\label{e:xaexp}\\
\vec v_A(t)&=&\vec v_{CM}(t)+ \vec\Omega(t)\times\Delta\vec x_A(t)
\nonumber\\
&&\qquad
+\sum_\mu \frac{d{\mathcal A}_\mu(t)}{dt}\,\mathbf{R}(t)\cdot \vec e^{\,\mu}_A, 
\qquad
\label{e:dxadtexp}
\end{eqnarray}
where $\Delta \vec x_A$ is given by
\begin{equation}
  \Delta \vec x_A(t) = {\mathbf{R}}(t)\cdot \left[\vec{x}_{0A}
    +\sum_\mu {\mathcal A}_\mu(t)\, \vec e^{\,\mu}_A \right].
\label{e:DeltaXDef}
\end{equation}

\subsection{Evolution Equations for Mode-Basis MD}
\label{s:md_large-scale}

The evolution equations for the mode-basis dynamical variables,
$\{\vec x_{CM},\mathbf{R}, \vec \Omega, \mathcal{A}_\mu\}$, are
determined from Eq.~(\ref{e:NewtonsLaw}) using the transformation given
in Eq.~(\ref{e:ModeBasisTransformation1}).  The resulting equations
can be written in the form:
\begin{eqnarray}
\frac{d^{\,2} \vec x_{CM}}{dt^2}
&=&-\sum_A\frac{1}{M}\frac{\partial U}{\partial \vec x_A},
\label{e:dxcmdt}\\
\frac{d\,\mathbf{R}}{dt}&=&-\mathbf{*\Omega}\cdot\mathbf{R},
\label{e:drdt}\\
\frac{d\,\vec\Omega}{dt}&=&
\tilde{\mathcal{U}}^{-1}\cdot\vec{\mathcal{V}}
,\label{e:dOmegadtEq}\\
\frac{d^{\,2}\!\mathcal{A}_\mu}{dt^2}&=&\sum_\nu {\mathcal S}^{\mu\nu}
\frac{d{\mathcal A}_\nu}{dt} + \sum_\nu {\mathcal T}^{\mu\nu}{\mathcal A}_\nu
+{\mathcal F}^\mu.\label{e:dAmuEq}
\end{eqnarray}
The tensor $\tilde{\mathcal{U}}$ and vector $\vec{\mathcal{V}}$ that
appear in Eq.~(\ref{e:dOmegadtEq}) are functions of $\{\vec
x_{CM},\mathbf{R}, \vec \Omega, \mathcal{A}_\mu,
d\mathcal{A}_\mu/dt\}$ given by
\begin{eqnarray}
\tilde{\mathcal{U}} &=& \sum_A \left(\frac{m_A}{M}\right)
\Bigl[\mathbf{R}^{-1}\, \bigl(\vec x_{0A}+\delta \vec x_A\bigr)\cdot\vec x_{0A}
  \nonumber\\
  &&\quad\qquad\qquad
  -\bigl(\vec x_{0A}+\delta \vec x_A\bigr)
  \otimes\bigl(\mathbf{R}\cdot\vec x_{0A}\bigr)\Bigr],
\label{e:UDef}\\
\vec{\mathcal{ V}} &=&\sum_A\left(\frac{m_A}{M}\right)\biggl[
2\Bigl(\vec x_{0A}\cdot \mathbf{R}^{-1}\cdot\vec\Omega\Bigr)
\frac{d\delta\vec x_A}{dt}\nonumber\\
&&\quad\qquad\qquad
-2\left(\vec x_{0A}\cdot\frac{d\delta\vec x_A}{dt}\right)
\left(\mathbf{R}^{-1}\cdot \vec \Omega\right)\nonumber\\
&&\quad\qquad\qquad
-\left(\vec\Omega\cdot\Delta \vec X_A\right) \vec x_{0A}\times
\left(\mathbf{R}^{-1}\cdot\vec\Omega\right)\nonumber\\
&&\quad\qquad\qquad
-\frac{1}{m_A}\vec x_{0A}\times\left(\mathbf{R}^{-1}\cdot
\frac{\partial U}{\partial \vec x_A}\right)
\biggr],\label{e:VDef}
\end{eqnarray}
where $\delta \vec x_A$ is given in Eq.~(\ref{e:ModeSumd}).  Similarly
the quantities $\mathcal{S}^{\mu\nu}$, $\mathcal{T}^{\mu\nu}$, and
$\mathcal{F}^{\mu\nu}$ that appear in Eq.~(\ref{e:dAmuEq}) are
functions of $\{\vec x_{CM},\mathbf{R}, \vec \Omega,
\mathcal{A}_\mu,d\mathcal{A}_\mu/dt\}$ given by
\begin{eqnarray}
\mathcal{S}^{\mu\nu}&=&2\sum_A \frac{m_A}{M}
\left(\mathbf{R}\cdot\vec e^{\,\mu}_A\right)\cdot
\left[\left(\mathbf{R}\cdot\vec e^{\,\nu}_A\right)\times\vec \Omega\right],
\label{e:SmunuDef}\\
\mathcal{T}^{\mu\nu}&=&\vec \Omega\cdot\vec \Omega\,\, \delta^{\mu\nu}
+\sum_A \frac{m_A}{M}
\left(\mathbf{R}\cdot\vec e^{\,\mu}_A\right)\cdot
\left[\left(\mathbf{R}\cdot\vec e^{\,\nu}_A\right)\times\frac{d\vec \Omega}{dt}
\right]\nonumber\\
&&
-\sum_A \frac{m_A}{M}
\left[\left(\mathbf{R}\cdot\vec e^{\,\mu}_A\right)\cdot\vec \Omega\right]
\left[\left(\mathbf{R}\cdot\vec e^{\,\nu}_A\right)\cdot\vec \Omega\right],
\label{e:TmunuDef}\\
\mathcal{F}^{\,\mu}&=&-\sum_A\frac{1}{M}\left(\mathbf{R}\cdot \vec e^{\,\mu}_A
\right)\cdot\frac{\partial U(\vec x_C)}{\partial \vec x_A}\nonumber\\
&&
+\sum_A\frac{m_A}{M}\left(\mathbf{R}\cdot\vec e^{\,\mu}_A\right)\cdot
\left[\left(\mathbf{R}\cdot\vec x_{0A}\right)\times\frac{d\vec \Omega}{dt}
\right]\nonumber\\
&&-\sum_A\frac{m_A}{M}\left\{\left[\left(\mathbf{R}\cdot\vec e^{\,\mu}_A\right)
\cdot\vec \Omega\right]
\left[\left(\mathbf{R}\cdot\vec x_{0A}\right)\cdot\vec \Omega\right]\right.
\nonumber\\
&&\qquad\qquad\qquad\qquad\;\;\left.
-\left(\vec e^{\,\mu}_A\cdot\vec x_{0A}\right)\left(\vec \Omega\cdot 
\vec \Omega\right)
\right\}.\qquad
\label{e:FmuDef}
\end{eqnarray}
The $d\,\vec\Omega/dt$ that appear in Eqs.~(\ref{e:TmunuDef}) and
(\ref{e:FmuDef}) are to be replaced by the expression on the right
side of Eq.~(\ref{e:dOmegadtEq}).  With those replacements, the
expressions on the right sides of
Eq.~(\ref{e:dxcmdt})--(\ref{e:dAmuEq}) depend only on $\{\vec
x_{CM},d\vec x_{CM}/dt, \mathbf{R}, \vec \Omega, \mathcal{A}_\mu,
d\mathcal{A}_\mu/dt\}$.

While the derivations leading to
Eqs.~(\ref{e:dxcmdt})--(\ref{e:dAmuEq}) are straightforward, they are
lengthy and have not been reproduced here in detail.  Those
derivations can be summarized however as follows.
Equation~(\ref{e:dxcmdt}) is obtained by inserting
Eq.~(\ref{e:NewtonsLaw}) into the second time derivative of
Eq.~(\ref{e:CMDef}).  Equation~(\ref{e:drdt}) follows trivially from
Eqs.~(\ref{e:AngularVelocityMatrix}) and (\ref{e:AngularVelocityDef}).
The derivation of Eq.~(\ref{e:dOmegadtEq}) is more complicated.  It is
obtained by projecting Eq.~(\ref{e:NewtonsLaw}) onto the three
independent generators of rigid rotations of the molecule, i.e.  the
vectors $\vec\theta\times \vec x_{0A}$ where $\vec\theta$ is a unit
vector whose direction determines the axis of rotation.  Similarly
Eq.~(\ref{e:dAmuEq}) is obtained by projecting
Eq.~(\ref{e:NewtonsLaw}) onto each of the mode-basis vectors $\vec
e^{\,\mu}_A$.

Finally, we note that while Eq.~(\ref{e:drdt}) determines the
evolution of the rotation matrix $\mathbf{R}(t)$, solving this
equation numerically directly in this form is problematic.  Instead it
is better to adopt some parameterization for the rotation matrices,
e.g. using Euler angles, and then to solve numerically the
differential equations implied by Eq.~(\ref{e:drdt}) for those
parameters.  In our numerical work we have adopted the ``quaternion''
parameterization of rotation matrices, in which each rotation matrix
is represented by four parameters $\{q_0,q_1,q_2,q_3\}$ with
$q_0^2+q_1^2+q_2^2+q_3^2=1$.
Appendix~\ref{s:QuaternionRepresentation} describes this quaternion
representation, and explicitly gives the representation of
Eq.~(\ref{e:drdt}) in terms of these parameters.  We point out that
the version of the quaternion evolution equations used here introduces
a new (so far as we know) constraint damping mechanism that ensures
the constraint, $q_0^2+q_1^2+q_2^2+q_3^2=1$, remains satisfied by the
numerical solution.

\subsection{Large-Scale MD Approximations}
\label{s:LargeScaleMD}

Equations~(\ref{e:dxcmdt})--(\ref{e:dAmuEq}) are a well-posed system
of ordinary differential equations that give an exact representation
of classical MD in terms of the mode-basis dynamical variables $\{\vec
x_{CM},\mathbf{R}, \vec \Omega, \mathcal{A}_\mu\}$.  The idea of our
large-scale MD approximations is to use some subset of the exact
equations to determine the macroscopic degrees of freedom, $\{\vec
x_{CM},\mathbf{R}, \vec \Omega\}$, while using simpler approximate
equations to determine the evolution of the internal vibrational
degrees of freedom $\mathcal{A}_\mu$.

The most straightforward way to construct a large-scale approximation
uses Eqs.~(\ref{e:dxcmdt})--(\ref{e:dOmegadtEq}) to determine $\{\vec
x_{CM},\mathbf{R}, \vec \Omega\}$, while replacing
Eq.~(\ref{e:dAmuEq}) with some approximate equation for
$\mathcal{A}_\mu$.  Perhaps the most natural approximation for
$\mathcal{A}_\mu$, which we refer to as the sinusoidal mode amplitude
(SMA) approximation, would be to set the mode amplitudes
$\mathcal{A}_\mu(t)$ to their small-amplitude perturbation solution
values:
\begin{eqnarray}
  \mathcal{A}_\mu(t) = \mathcal{A}^0_\mu\sin(\omega_\mu t+ \varphi_\mu),
  \label{e:sinusoidal_approx}             
\end{eqnarray}
where $\omega_\mu$ is the mode frequency determined from
Eq.~(\ref{e:ModeEquation}), while $\mathcal{A}_\mu^0$ and
$\varphi_\mu$ are constants that specify the amplitude of phase of
each mode.  In this approximation, Eq.~(\ref{e:sinusoidal_approx})
replaces Eq.~(\ref{e:dAmuEq}) and is used to evaluate the right sides
of Eqs.~(\ref{e:dxcmdt})--(\ref{e:dOmegadtEq}).  Those equations for
the macroscopic degrees of freedom $\{\vec x_{CM},\mathbf{R}, \vec
\Omega\}$ are then solved numerically.  The system of equations being
solved numerically is therefore reduced from $6N$ first-order
equations for the exact MD system, to just twelve for this large-scale
approximation.  The use of this approximation eliminates the need to
evaluate the complicated quantities $\mathcal{S}^{\mu\nu}$,
$\mathcal{T}^{\mu\nu}$ and $\mathcal{F}^{\mu}$ that appear on the
right side of Eq.~(\ref{e:dAmuEq}) numerically.  This reduction in the
number of equations to be solved numerically, as well as the reduction
in the need to evaluate the complicated expressions that occur on the
right side of Eq.~(\ref{e:dAmuEq}) considerably reduces the
computational cost of implementing the SMA approximation.

Another plausible approximation, which we refer to as the zero mode
amplitude (ZMA) approximation, simply sets all the mode amplitudes to
zero:
\begin{eqnarray}
  \mathcal{A}_\mu(t) = 0.
  \label{e:zero_approx}             
\end{eqnarray}
We expect the long time averages of the positions of the atoms to be
their equilibrium positions $\vec x_{0A}$.  Thus we expect the time
averages of the mode amplitudes $\mathcal{A}_\mu(t)$ to be zero, as
they are for example in the SMA approximation given in
Eq.~(\ref{e:sinusoidal_approx}).  This should reduce the computational
cost of evaluating the right sides of
Eqs.~(\ref{e:dxcmdt})--(\ref{e:dOmegadtEq}) in the ZMA approximation
even below those costs in the SMA approximation.  In addition the ZMA
approximation eliminates all the short timescale effects associated
with the molecular vibrations, so it should be possible to use much
larger timesteps to determine the macroscopic variables $\{\vec
x_{CM},\mathbf{R}, \vec \Omega\}$ in this approximation, and thus
to reduce the computational cost even below those for the SMA
approximation.

We note that hybrid approximations can easily be constructed
as well.  In these hybrid approximations some of the mode amplitudes
are set to the SMA or the ZMA approximations given in
Eqs.~(\ref{e:sinusoidal_approx}) or (\ref{e:zero_approx}), while the
remaining amplitudes are determined numerically by solving
Eq.~(\ref{e:dAmuEq}).  This approach might be appropriate for systems
having a few modes with oscillation timescales comparable to the
timescales associated with the macroscopic properties of the molecule.
In such cases those low frequency modes could be treated exactly while
the approximations could still be used for the majority of modes
having much shorter oscillation timescales.

Somewhat more sophisticated large-scale approximations can also be
obtained by choosing the equations of motion for the macroscopic
variables, $\{\vec x_{CM},\mathbf{R}, \vec \Omega\}$, from the
particular combination of the exact equations that determine the
evolution of the macroscopic linear and angular momentum, $\vec P$ and
$\vec J$, of each molecule.  These quantities are defined by
\begin{eqnarray}
  \vec P &=& \sum_A m_A\frac{d\vec x_A}{dt},
  \label{e:PDef}\\
  \vec J &=& \sum_A m_A \Delta \vec x_A\times\frac{d\Delta\vec x_A}{dt},
  \label{e:JDef}
\end{eqnarray}
where $\Delta \vec x_A$, defined in Eq.~(\ref{e:xaexp}), is the
position of each atom relative to the center of mass of the
molecule.
\footnote{The angular momentum $\vec J$ defined in this way
  is the angular momentum relative to the center of mass of the
  molecule.}
The time derivatives of these quantities can be written
in terms of the macroscopic variables as
\begin{eqnarray}
  \frac{d\vec P}{dt}&=& M\frac{d^{\,2}\vec x_{CM}}{dt^2},
  \label{e:Pdot}\\
  \frac{d\vec J}{dt}&=& 
  \tilde{\mathcal{J}}_\Omega\cdot\frac{d\vec\Omega}{dt}+\vec{\mathcal{J}},
  \label{e:Jdot}
\end{eqnarray}
where  $\tilde{\mathcal{J}}_\Omega$ and $\vec{\mathcal{J}}$
are given by
\begin{eqnarray}
  \tilde{\mathcal{J}}_\Omega &=& \sum_A m_A
  \Bigl(\Delta \vec x_A\cdot\Delta\vec x_A
  \,\mathbf{I}- \Delta \vec x_A\otimes \Delta\vec x_A\Bigr),
  \label{e:JmatrixDef}\\
  \vec{\mathcal{J}}&=&\sum_A m_A \Delta\vec x_A \times \vec{\mathcal{B}}_A,
  \label{e:JvectorDef}
\end{eqnarray}
and where $\vec{\mathcal{B}}_A$ is given by
\begin{eqnarray}
  \vec{\mathcal{B}}_A &=& 
  \sum_\mu\left(\frac{d^{\,2}\mathcal{A}_\mu}{dt^2}
  +2\frac{d\mathcal{A}_\mu}{dt}\,
  \vec\Omega\times\right)\mathbf{R}\cdot\vec e^{\,\mu}_A
    \nonumber\\
    &&
    \qquad
  +\left(\vec\Omega\otimes\vec\Omega -\vec\Omega\cdot\vec\Omega\,\mathbf{I}
  \right)\cdot\Delta \vec x_A.\qquad
\end{eqnarray}

If we assume the mode amplitudes $\mathcal{A}_\mu(t)$ are
predetermined by some approximate expressions, like those in
Eqs.~(\ref{e:sinusoidal_approx}) or (\ref{e:zero_approx}) for example,
then $\tilde{\mathcal{J}}_\Omega$ and $\vec{\mathcal{J}}$ depend only
on the large scale variables $d\vec x_{CM}/dt$, $\mathbf{R}$,
$\vec\Omega$, but not on their time derivatives.
Equations~(\ref{e:Pdot}) and (\ref{e:Jdot}) can therefore be used to
construct an alternate set of approximate evolution equations for the
large scale variables.  In particular an evolution equation for the
center-of-mass motion of each molecule can be obtained by setting the
rate of change of the total momentum equal to the total external force
acting on the molecule:
\begin{eqnarray}
  \frac{d\vec P}{dt}&=&-\sum_A\frac{\partial U}{\partial \vec x_A}
  = M\frac{d^{\,2}\vec x_{CM}}{dt^2}.
\end{eqnarray}
Similarly an equation for $d\vec\Omega/dt$ can be obtained by setting
$d\vec J/dt$ equal to the total external torque acting on the
molecule:
\begin{eqnarray}
  \frac{d\vec J}{dt}&=& -\sum_A \Delta\vec x_A \times
  \frac{\partial U}{\partial \vec x_A}=
  \tilde{\mathcal{J}}_\Omega\cdot\frac{d\,\vec\Omega}{dt}+\vec{\mathcal{J}}.
  \label{e:AngMomCons}
\end{eqnarray}
The resulting evolution equations for $\vec x_{CM}$ and
$\vec\Omega$ are given by
\begin{eqnarray}
  \frac{d^{\,2}\vec x_{CM}}{dt^2}&=& -\frac{1}{M}\sum_A\frac{\partial
    U}{\partial \vec x_A},
  \label{e:NewXcmEq}\\
     \frac{d\,\vec\Omega}{dt}
   &=&-\Bigl(\tilde{\mathcal{J}}_\Omega \Bigr)^{-1}\cdot\left[
   \vec{\mathcal{J}} +
   \sum_A \Delta\vec x_A \times
   \frac{\partial U}{\partial \vec x_A}\right].\qquad
     \label{e:NewOmegaEq}
\end{eqnarray}
Equations~(\ref{e:NewXcmEq}) and (\ref{e:NewOmegaEq}) represent
somewhat different projections of the exact MD equations,
Eq.~(\ref{e:NewtonsLaw}), than those given in Eqs.~(\ref{e:dxcmdt})
and (\ref{e:dOmegadtEq}).  Therefore Eqs.~(\ref{e:NewXcmEq}) and
(\ref{e:NewOmegaEq}) together with Eq.~(\ref{e:drdt}), provide an
alternate somewhat different set of evolution equations for the
macroscopic variables $\{\vec x_{CM},\mathbf{R}, \vec \Omega\}$.  We
refer to these alternate equations as the momentum conserving (MC)
large-scale approximation.  These equations can be solved using any
predetermined approximate form for the mode amplitudes
$\mathcal{A}_\mu(t)$.  In this paper we explore the two possibilities
discussed above: We refer to the momentum conserving approximation
using sinusoidal mode amplitudes, Eq.~(\ref{e:sinusoidal_approx}), as
the MCSMA approximation, and the momentum conserving approximation
using zero mode amplitudes, Eq.~(\ref{e:zero_approx}), as the MCZMA
approximation.

\section{Reliability Testing}
\label{s:ReliabilityAssessment}

In this section we assess the reliability of the large-scale MD
approximations introduced in Sec.~\ref{s:LargeScaleMD}.  We do this by
comparing the values of the macroscopic variables, $\{\vec x_{CM},\vec
v_{CM}, \mathbf{R}, \vec\Omega\}$, computed numerically using the
exact MD equations, with those computed using several examples of
large-scale MD approximations.  We also compare how well these various
methods conserve the total energy $E$, the total momentum $\vec P$,
and the total angular momentum $\vec J$ of each molecule.  The
remainder of this section is organized as follows.
Section~\ref{s:ModelProblem} describes in detail the model problem,
and the numerical methods used to solve the MD equations for these
tests.  Section~\ref{s:TestingModeBasisRepresentation} presents the
results of numerical tests that confirm the mode-basis representation
of the exact MD equations, introduced in Sec.~\ref{s:md_large-scale},
gives the same results as the standard Cartesian representation for
this model problem.  Section~\ref{s:TestingLargeScaleApproximations}
gives the results of our numerical tests for the large-scale
approximations SMA, ZMA, MCSMA and MCZMA introduced in
Sec.~\ref{s:LargeScaleMD}.  Finally in
Sec.~\ref{s:ComputationalEfficiency} we compare the computational
efficiency of these various methods when applied to our model problem.

 \subsection{Model Problem}
 \label{s:ModelProblem}

\begin{figure}[ht]
\includegraphics[width=3.5cm]{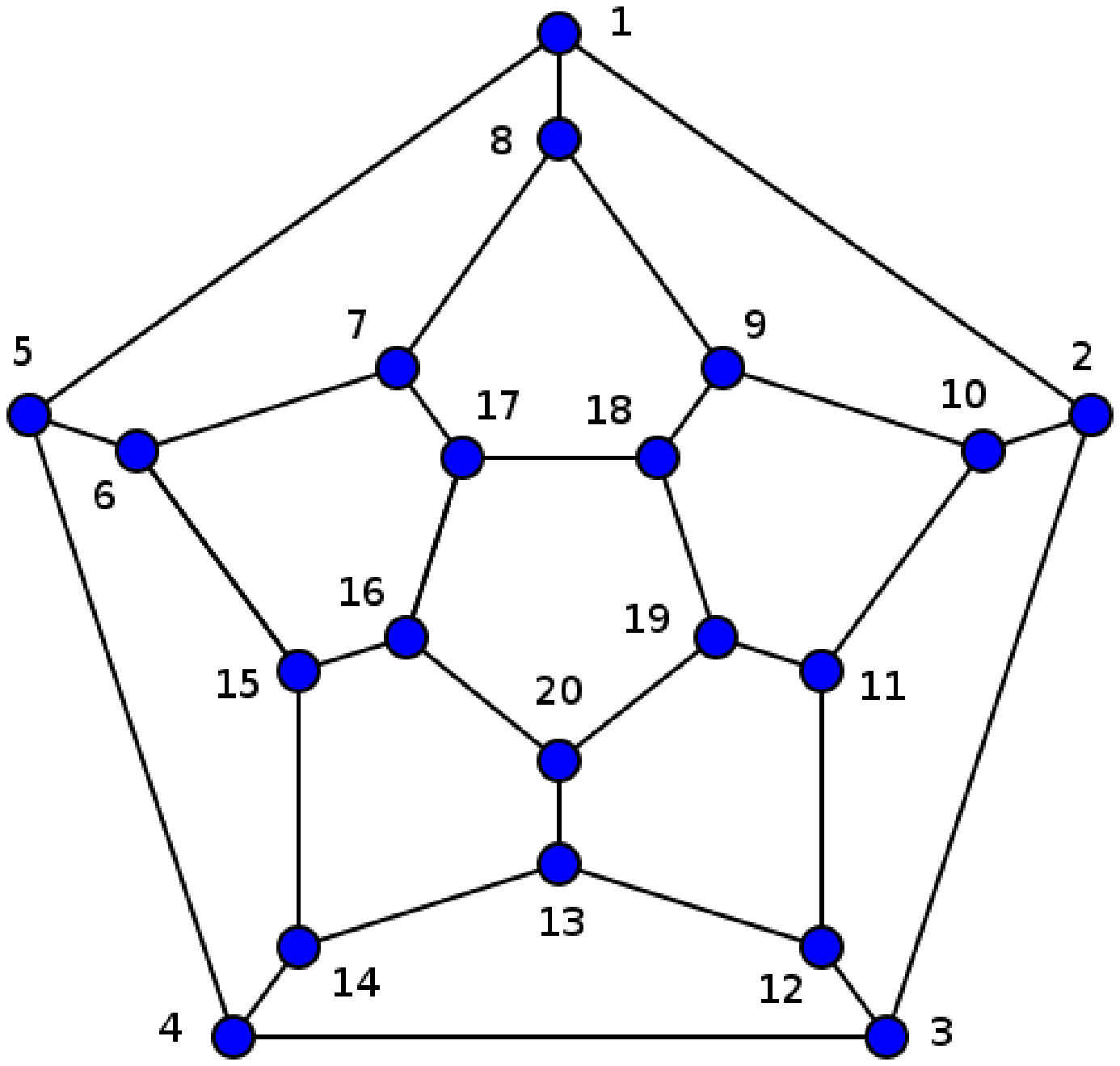}
\hspace{0.8cm}
\includegraphics[width=3.5cm]{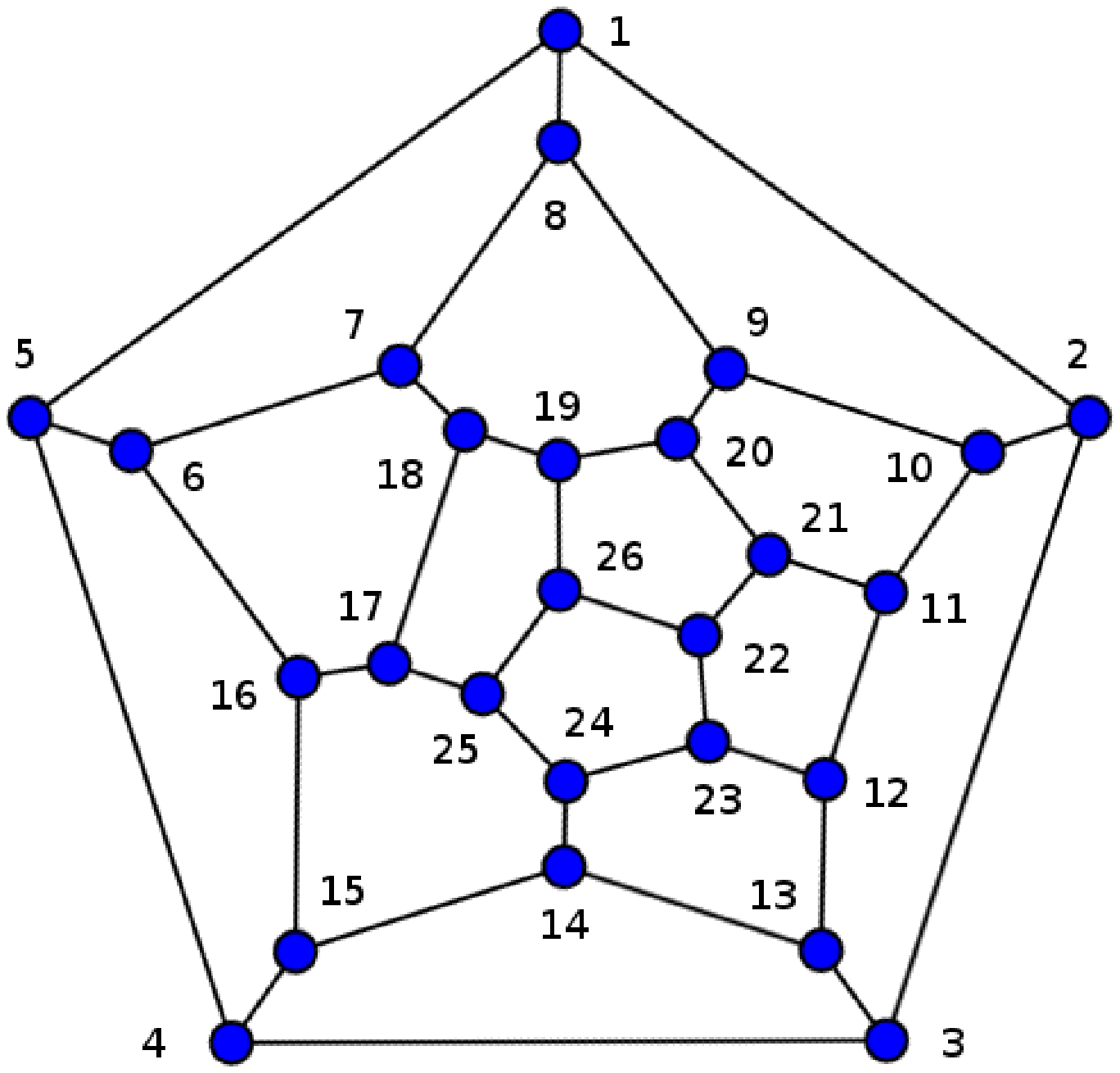}
\hspace{0.8cm}
\includegraphics[width=3.5cm]{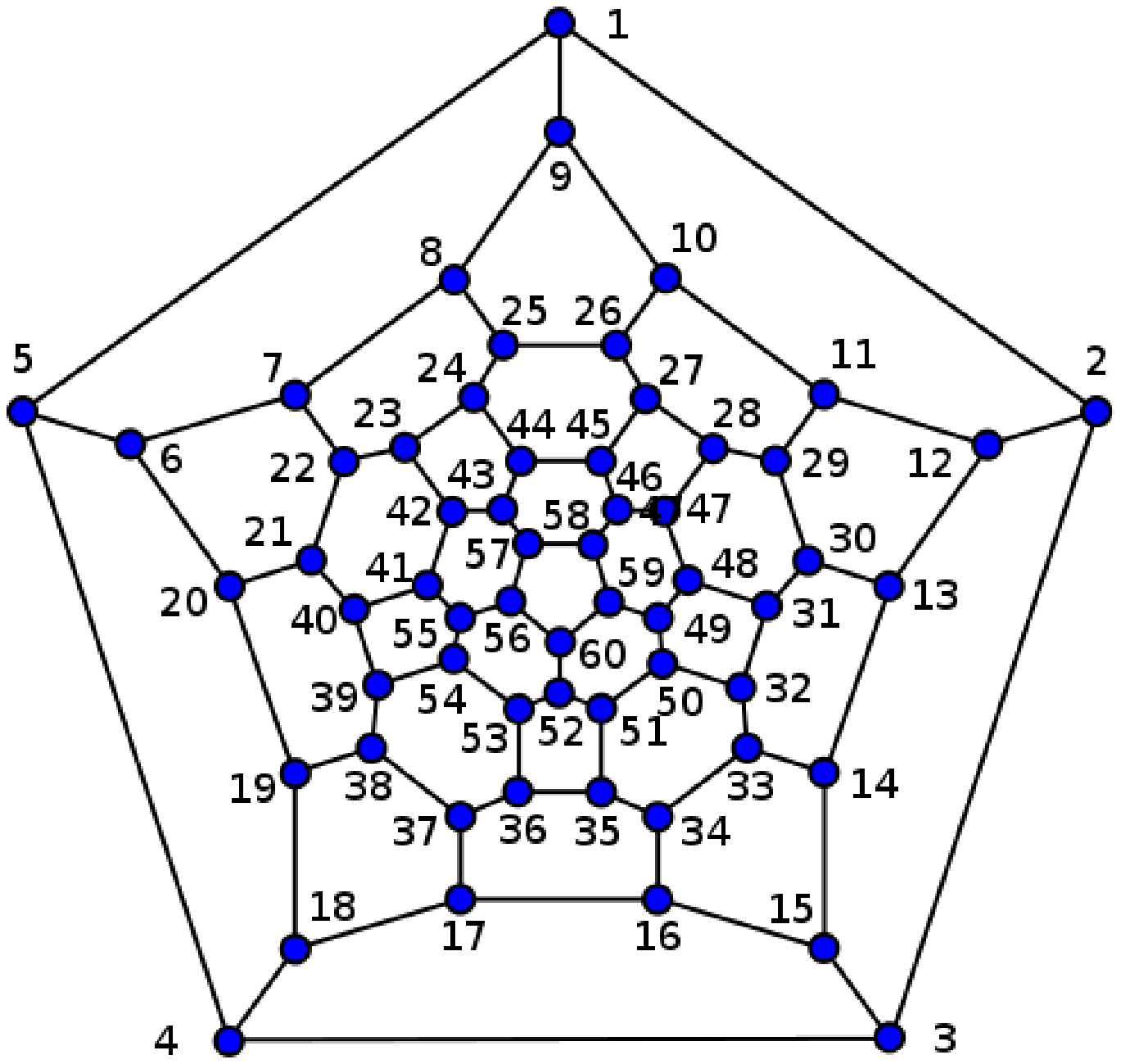}
\hspace{0.8cm}
\includegraphics[width=3.5cm]{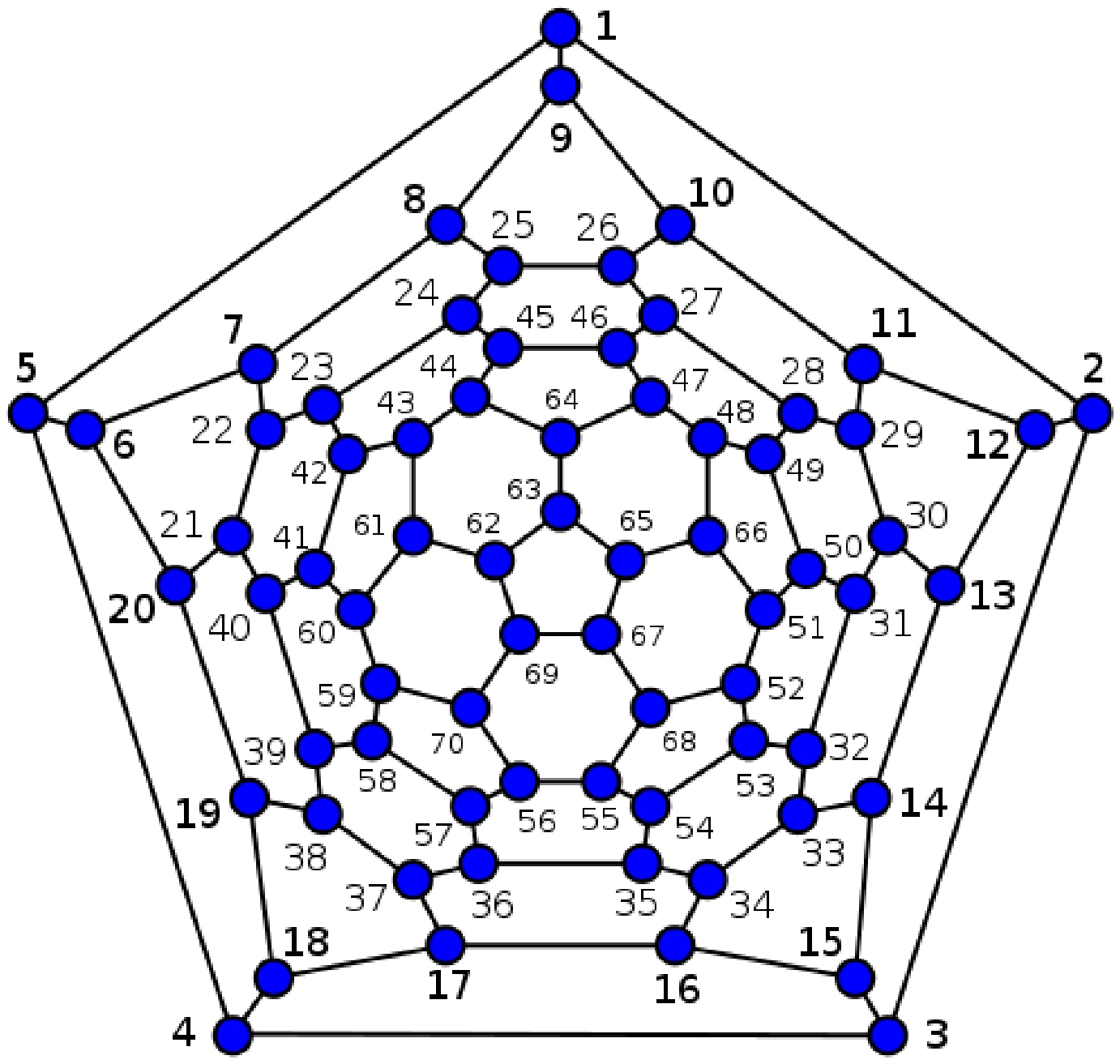}
\caption{Figures show the topological bond connections for the
  $C_{20}$, $C_{26}$, $C_{60}$, and $C_{70}$ molecules used in our
  test MD simulations.  Each atom is labeled with a number that
  represents the value of the index $\scriptstyle A$ for that atom.}
\label{f:FullerineTopology}
\end{figure}

We use a collection of simple molecules, the fullerenes $C_{20}$,
$C_{26}$, $C_{60}$, and $C_{70}$, to study the reliability of the
large-scale MD approximations introduced in Sec.~\ref{s:LargeScaleMD}.
These molecules consist entirely of trivalent carbon atoms located at
the vertices of convex polyhedra.  Figure~\ref{f:FullerineTopology}
illustrates the topological bond connections between the atoms (but
not the actual geometrical shapes) of the molecules used in our tests.
 
We use a simplified version of the CHARMM27 model for the potential
energy function $U(\vec x_A)$ that determines the interactions between
the atoms.  This interaction potential is given by
\begin{eqnarray}
U&=&\tfrac{1}{2}\kappa_b\!\!\!\!\!\!
\sum_{\mathrm{bonds}\,\,(CD)}\!\!\!\!(r_{CD}-L_b)^2
+\tfrac{1}{2}\kappa_\theta\!\!\!\!\!\!\!\!\! 
\sum_{\mathrm{angles}\,\,(CDE)}\!\!\!\!\!\! (\theta_{CDE}-\theta_b)^2,
\nonumber\\\label{e:energydef}
\end{eqnarray}
where $L_b$ and $\theta_b$ are the lengths and angles of equilibrium
molecular bonds in this simple model, and where $r_{CD}$ and
$\theta_{CDE}$ represent respectively the distance between atoms
$\scriptstyle C$ and $\scriptstyle D$ and the angle formed by the bonds between
atom $\scriptstyle D$ with atoms $\scriptstyle C$ and $\scriptstyle E$,
\begin{eqnarray}
r_{CD}^2 &=& (\vec x_C - \vec x_D)
\cdot(\vec x_C - \vec x_D),\\
\cos\theta_{CDE}&=&\frac{(\vec x_C-\vec x_D)\cdot(\vec x_E - \vec x_D)}
{r_{CD}\,r_{DE}}.
\end{eqnarray}
The sums in Eq.~(\ref{e:energydef}) are taken respectively over the
collection of bonds $\scriptstyle(CD)$ between the pairs of atoms
$\scriptstyle C$ and $\scriptstyle D$, and over the collection of
angles $\scriptstyle(CDE)$ formed by the bond between atoms
$\scriptstyle C$ and $\scriptstyle D$ and the bond between atoms
$\scriptstyle D$ and $\scriptstyle E$.  We use the following values
for the carbon-carbon bond parameters: $\kappa_b = 305$ kcal/
\hspace{-1.5mm}\AA$^2$/mol, $L_b = 1.375$ \AA, and $ \theta_b = 120$
degrees taken from the CHARMM27 force field parameters for these
carbon-carbon bonds.~\cite{MacKerell1998} For simplicity in coding up
our tests, we left out the standard torsion-angle bond interactions
from the potential energy.  Without those torsion-angle forces, the
fullerene molecules are unstable with our simplified interaction
potential when the CHARMM27 value is used for the bond angle force
constant, $\kappa_\theta = 40$ kcal/rad$^2$/mol.  Therefore we have
increased the value of $\kappa_\theta$ used in our tests to
$\kappa_\theta = 305$ kcal/rad$^2$/mol to achieve stability.  Our
purpose here is to test the robustness of our large-scale
approximations.  These approximations should succeed or fail
independent of the details of the interaction potential model being
used, so we do not think it matters that our simplified potential
model is not state of the art.

The first step in our analysis of these molecules is to determine
their equilibrium configurations, $\vec x_{0A}$, for the interaction
potential $U(\vec x_A)$ given in Eq.~(\ref{e:energydef}).  We do this
by finding the energy minimum where $\partial U/\partial\vec x_A=0$.
We use the Fletcher-Reeves-Polak-Ribiere version of the conjugate
gradient method with line minimizations to find these equilibrium
states, $\vec x_{0A}$, numerically.~\cite{numrec_f} Given an
equilibrium state, we next evaluate the Hessian matrix $\partial^2
U/\partial\vec x_A\partial \vec x_B$ numerically for that state, and
solve Eq.~(\ref{e:ModeEquation}) to determine the eigenvalues
$\omega_\mu$ and eigenvectors $\vec e^{\,\mu}_A$.  We use Householder
reduction to transform $\partial^2 U/\partial\vec x_A\partial \vec
x_B$ to tridiagonal form, followed by a traditional QL algorithm with
implicit shifts to determine the eigenvalues and eigenvectors
numerically.~\cite{numrec_f} These eigenvectors are then projected and
normalized so they satisfy
Eqs.~(\ref{e:ModeOrthogonality})--(\ref{e:RotationConstraint}) to
double precision accuracy numerically.

We construct initial data for our test evolutions by choosing values
for the mode-basis variables $\{\vec x_{CM},\vec v_{CM}, \mathbf{R},
\vec\Omega,\mathcal{A}_\mu,d\mathcal{A}_\mu/dt\}$ that are appropriate
for a thermodynamic equilibrium state at temperature $T$.  Following
the equipartition theorem, we fix the values for each of the
mode-basis variables so that each degree of freedom of the molecule
has energy $\half kT$, where $k = 1.9872\times 10^{-3}$ kcal/(mol K)
is Boltzmann's constant.  All the tests reported here use a
temperature $T=300K$.  By choosing the origin and the orientation of
the Cartesian coordinate system, we can set $\vec x_{CM}=0$ and
$\mathbf{R}=\mathbf{I}$ at $t=0$ without loss of generality.  We
choose $\vec v_{CM}$ and $\vec \Omega$ at $t=0$ to be vectors whose
orientations are set with a random number generator, and whose
magnitudes are set by requiring the translational and rotational
kinetic energies to satisfy,
\begin{eqnarray}
  \half kT &=& \half M\, \vec v_{CM}\cdot\vec v_{CM},\\
  \half kT &=& \half \sum_A m_A \,\Bigl[\vec x_{0A}\cdot\vec x_{0A}
    \,\,\vec\Omega\cdot\vec\Omega-
    \bigl(\vec x_{0A}\cdot\vec\Omega\bigr)^2\Bigr].
\end{eqnarray}
The mode amplitudes $\mathcal{A}_\mu$ and their time derivatives
$d\mathcal{A}_\mu/dt$ are chosen at $t=0$ to ensure that each normal
mode of the molecule has energy $kT$:
\begin{eqnarray}
  \mathcal{A}_\mu&=& \frac{1}{\omega_\mu}\sqrt{\frac{2kT}{M}}\sin\varphi_\mu,\\
  \frac{d\mathcal{A}_\mu}{dt}&=&\sqrt{\frac{2kT}{M}}\cos\varphi_\mu,
\end{eqnarray}
where $\varphi_\mu$ are randomly selected phases.  These initial
values for $\{\vec x_{CM},\vec v_{CM}, \mathbf{R},
\vec\Omega,\mathcal{A}_\mu,d\mathcal{A}_\mu/dt\}$ are used to start
the evolutions of the exact mode-basis representation of the MD
equations.  We convert them to the equivalent Cartesian representation
variables using Eqs.~(\ref{e:xaexp}) and (\ref{e:dxadtexp}), and use
those initial data to start our exact Cartesian MD evolutions.  We
also use the same initial values of $\{\vec x_{CM},\vec v_{CM},
\mathbf{R}, \vec\Omega,\mathcal{A}_\mu,d\mathcal{A}_\mu/dt\}$ to set
the initial data for the SMA and MCSMA approximation tests.  And
finally we use these same values for $\{\vec x_{CM},\vec v_{CM},
\mathbf{R}, \vec\Omega\}$ with $\mathcal{A}_\mu=d\mathcal{A}_\mu/dt=0$
to set the initial data for the ZMA and MCZMA approximation tests.

All the MD simulation methods considered here consist of systems of
ordinary differential equations of the form $d\,\vec Y/dt = \vec
F(\vec Y,t)$, where $\vec Y$ is the n-dimensional vector consisting of
the dynamical fields in a particular method and $\vec F(\vec Y,t)$ is
the right side of the evolution equations for those fields.  We solve
these systems numerically with the initial data $\vec Y=\vec Y(0)$
described above using an 8th order integrator by Dormand and Prince
with dynamic timestep size control (see Hairer et
al.~\cite{Hairer2008} for details).  This algorithm controls the error
in $\vec Y(t)$ by adjusting the timestep size to keep an estimate of
the local time-truncation error below $\epsilon_\tau\bigl(|\vec
Y|+1\bigr)$ (see Hairer et al.~\cite{Hairer2008} for details about
this timestep control).  We run each simulation with several values of
the timestep accuracy parameter $\epsilon_\tau$ in the range
$10^{-13}\leq\epsilon_\tau\leq 10^{-6}$ to verify that timestep errors
are not the dominant cause of any differences we may see between the
various MD evolution methods.

\subsection{Testing the Exact Mode-Basis Representation}
\label{s:TestingModeBasisRepresentation}

Our first numerical tests of the model problem described in
Sec.~\ref{s:ModelProblem} are designed to examine the differences
between MD simulations performed with the standard Cartesian-basis
representation of the MD equations, Eq.~(\ref{e:NewtonsLaw}), and the
exact mode-basis representation given in
Eqs.~(\ref{e:dxcmdt})--(\ref{e:dAmuEq}). To perform these tests we use
the exact Cartesian-basis solution computed with timestep accuracy
parameter $\epsilon_\tau=10^{-13}$ as the reference solution with
which to compare everything else.  We refer to this reference
solution as $\vec x^{\,Ref}_A(t)$.

\begin{figure}[!t]
\includegraphics[width=.45\textwidth]{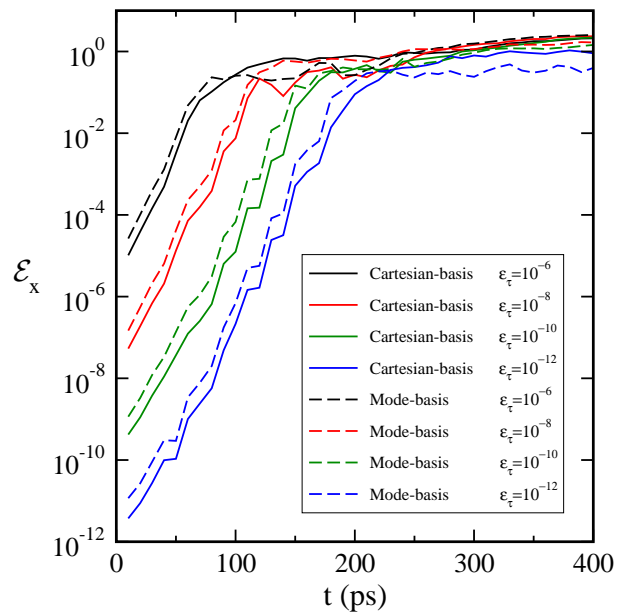}
\caption{Solid curves show $\mathcal{E}_x$ computed with the exact
  Cartesian-basis MD code for different values of the timestep
  accuracy parameter $\epsilon_\tau$.  Dashed curves show
  $\mathcal{E}_x$ for evolutions of the exact mode-basis MD code for
  different $\epsilon_\tau$. The exponential growth in these
  $\mathcal{E}_x$ is caused by the chaotic nature of MD evolutions.}
\label{f:XConvergence}
\end{figure}
We first measure how sensitively the exact Cartesian-basis MD
evolutions depend on the timestep accuracy parameter $\epsilon_\tau$.
To do this we evaluate the quantity $\mathcal{E}_x$ that measures the
differences between solutions computed using different timestep
accuracy parameters, $\epsilon_\tau$, and the reference solution:
\begin{eqnarray}
  \mathcal{E}_x(\epsilon_\tau)=\sqrt{\frac{1}{N}
    \sum_A\bigl|\vec x_A(\epsilon_\tau)-\vec x_A^{\,Ref}\bigr|^2}.
\end{eqnarray}
The solid curves in Fig.~\ref{f:XConvergence} show
$\mathcal{E}_x(\epsilon_\tau)$ for exact Cartesian-basis MD
simulations of the $C_{20}$ fullerene molecule computed with four
different values of
$\epsilon_\tau=\{10^{-6},10^{-8},10^{-10},10^{-12}\}$.  Data points
for these curves are obtained by evaluating $\mathcal{E}_x$ at 10 ps
time intervals during the evolutions.  Each of these curves begins at
early times with $\mathcal{E}_x\approx\epsilon_\tau$ and then grows
exponentially until $\mathcal{E}_x\approx 1$ where it remains
relatively constant for the duration of the simulation.  These curves
confirm the expectation that MD simulations are chaotic.  Although the
initial data at $t=0$ for these various runs are identical, after one
time step the solutions differ from the one specified by the initial
data, by amounts that depend on the timestep accuracy parameter
$\epsilon_\tau$.  By definition, chaotic dynamical systems have the
property that nearby solutions diverge exponentially.
Figure~\ref{f:XConvergence} confirms that this is what is going on by
showing that each of these evolutions of the fullerene $C_{20}$
molecule diverges from the reference solution at the same exponential
rate.  The analogous graphs for the other fullerene molecules,
$C_{26}$, $C_{60}$ and $C_{70}$, included in our study are very
similar, except for the timescale on which the chaotic instability
occurs.  In the $C_{26}$ case, the instability grows at about half the
rate of the $C_{20}$ case, while the instabilities in the $C_{60}$ and
the $C_{70}$ cases grow at two or three times the $C_{20}$ rate.  The
presence of chaotic behavior in these MD simulations demonstrates why
it is impossible to compute molecular evolutions in complete detail.
Only certain macroscopic features of the evolutions, like dynamically
conserved quantities such as the energy, momentum and angular momentum
are reproducible and simulatable.

The dashed curves in Fig.~\ref{f:XConvergence} show $\mathcal{E}_x$
for simulations based on the exact mode-basis MD representation.  We
use the same Cartesian-basis reference solution $\vec x_{A}^{\,Ref}$
when computing $\mathcal{E}_x$ for these mode-basis simulations, and
we see that the exponential divergence from the reference solution has
exactly the same structure it has for the Cartesian-basis evolutions.
The only difference is that the values of $\mathcal{E}_x$ are somewhat
larger, $\mathcal{E}_x\approx 3\epsilon_\tau$, at very early times in
the mode-basis case.  These differences appear to be caused by the
fact that the mode-basis equations are much more complicated than
their Cartesian-basis counterparts, so the truncation errors are
somewhat higher in this case for fixed $\epsilon_\tau$.  The rate of
the exponential divergence from the reference solution is the same as
the Cartesian-basis case, so Fig.~\ref{f:XConvergence} confirms that
the mode-basis equations produce the same evolutions as the standard
Cartesian-basis representation of MD.  More importantly perhaps, these
tests also confirm that our codes to evolve both versions of the MD
equations contain no serious errors.

\begin{figure}[!t]
  \includegraphics[width=.45\textwidth]{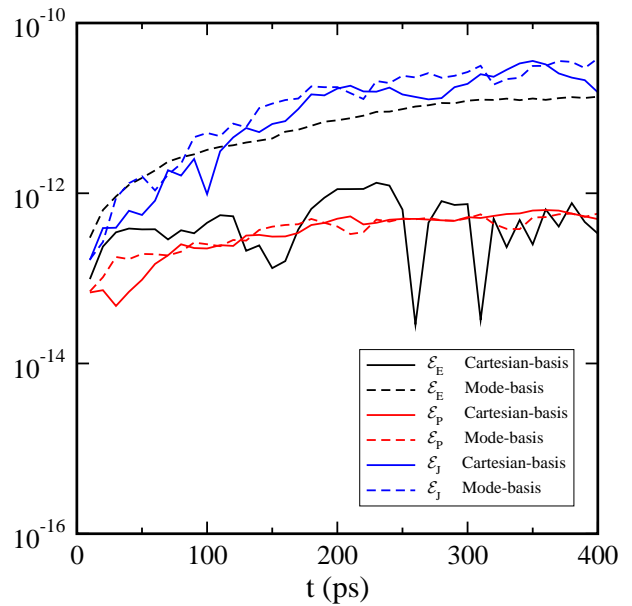}
  \caption{Solid curves show $\mathcal{E}_E$, $\mathcal{E}_P$ and
    $\mathcal{E}_J$ for evolutions using the exact Cartesian MD code
    with timestep accuracy $\epsilon_\tau=10^{-13}$.  Dashed curves
    show these same quantities for evolutions using the exact mode-basis MD
    code.}
\label{f:EnergyConservation}
\end{figure}

We have also monitored how well the energy $E$ defined by,
\begin{equation}
  E=\tfrac{1}{2}\sum_A \frac{d\vec x_A}{dt}\cdot\frac{d\vec x_A}{dt}+U(\vec x_B),
\end{equation}
the total momentum $\vec P$ defined in Eq.~(\ref{e:PDef}) and the
total angular momentum (about the center of mass) $\vec J$ defined in
Eq.~(\ref{e:JDef}) are conserved in these exact MD evolutions.  In the
absence of external forces (like van der Waals interactions with other
molecules) these quantities should all be conserved by the exact MD
evolution equations.  To monitor these conserved quantities we define
$\mathcal{E}_E$, $\mathcal{E}_P$ and $\mathcal{E}_J$ that measure the
deviations of these quantities from their initial values:
\begin{eqnarray}
  \mathcal{E}_E(t)&=&\frac{\bigl|E(t)-E(0)\bigr|}{E(0)},\label{e:Eerror}\\
  \mathcal{E}_P(t)&=&\frac{\bigl|\vec P(t)-\vec P(0)\bigr|}
    {\bigl|\vec P(0)\bigr|},\label{e:Perror}\\
  \mathcal{E}_J(t)&=&\frac{\bigl|\vec J(t)-\vec J(0)\bigr|}
    {\bigl|\vec J(0)\bigr|}.\label{e:Jerror}
\end{eqnarray} 
The solid curves in Fig.~\ref{f:EnergyConservation} show the evolution
of these energy and momentum conservation errors for the simulations
of the Cartesian-basis representation of the MD equations.  The
results in this figure represent those of the highest resolution
simulations, i.e. those computed with timestep error parameter
$\epsilon_\tau=10^{-13}$.  We see from these figures that each of the
conservation error quantities begins at small times with
$\mathcal{E}_E\approx \mathcal{E}_P\approx \mathcal{E}_J\approx
\epsilon_\tau$, which then grow slowly, roughly as a power law in
time: $\mathcal{E}\propto t^{\,k}$, with $k\lesssim 2$. (The growth of
truncation level errors in this way is typical of explicit numerical
ordinary differential equation integrators such as the Dormand-Prince
algorithm used in our tests.)  The dashed curves in
Fig.~\ref{f:EnergyConservation} show the errors in these conserved
quantities for evolutions of the same initial data using the mode-basis
representation of MD.  The mode-basis results shown in
Fig.~\ref{f:EnergyConservation} were also computed using timestep
error parameter $\epsilon_\tau=10^{-13}$.  These results confirm that
both versions of the exact MD equations conserve the energy, the
linear momentum, and the angular momentum of molecules at the level of
the numerical truncation error used.

\subsection{Testing the Large-Scale Approximations}
\label{s:TestingLargeScaleApproximations}

In this section we present the results of numerical tests of the
large-scale MD approximations developed in Sec.~\ref{s:LargeScaleMD}.
These approximations include the SMA and ZMA approximations that use a
sinusoidal-in-time approximation (SMA) or the zero approximation (ZMA)
respectively for the mode amplitudes $\mathcal{A}_\mu(t)$.  These
approximations solve the exact mode-basis evolution
Eqs.~(\ref{e:dxcmdt})--(\ref{e:dOmegadtEq}) for $\{\vec x_{CM}, \vec
v_{CM}, \mathbf{R},\vec\Omega\}$, and simply ignore the exact
evolution Eq.~(\ref{e:dAmuEq}) for $\mathcal{A}_\mu(t)$.  We also test
approximations that use angular momentum conservation instead of
Eq.~(\ref{e:dOmegadtEq}) to determine the evolution of $\vec
\Omega(t)$.  The equations for these momentum conserving
approximations, MCSMA and MCZMA, are given in Eqs.~(\ref{e:NewXcmEq})
and (\ref{e:NewOmegaEq}).  All the numerical results shown here use
the highest time resolution, $\epsilon_\tau=10^{-13}$, in evolutions
of the $C_{20}$ fullerene molecule.

\begin{figure}[!t]
  \includegraphics[width=.45\textwidth]{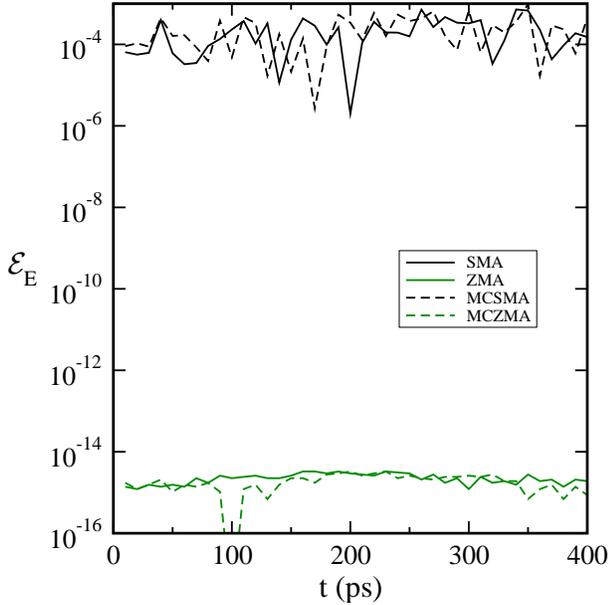}
  \caption{Solid curves show energy conservation violations
    $\mathcal{E}_E$ for the large-scale approximations SMA and ZMA,
    while dashed curves show these violations for the MCSMA and MCZMA
    approximations.}
\label{f:XMA_Energy}
\end{figure}
First we test how well these large-scale approximations conserve the
energy $E$, linear momentum $\vec P$, and angular momentum $\vec J$
during the evolutions of our model problem.  We use the quantities
$\mathcal{E}_E$, $\mathcal{E}_P$ and $\mathcal{E}_J$ defined in
Eqs.~(\ref{e:Eerror})--(\ref{e:Jerror}) to monitor conservation
violations.  The solid curves in Fig.~\ref{f:XMA_Energy} show the
energy conservation violations $\mathcal{E}_E$ for the SMA and ZMA
approximations, while the dashed curves show these violations for the
MCSMA and MCZMA approximations.  We see that the zero mode
approximations ZMA and MCZMA conserve energy much better than the
sinusoidal mode approximations SMA and MCSMA.  However, even these
sinusoidal mode approximations give energy conservation violations
below the $0.1\%$ level for these test problems.

\begin{figure}[!h]
  \includegraphics[width=.45\textwidth]{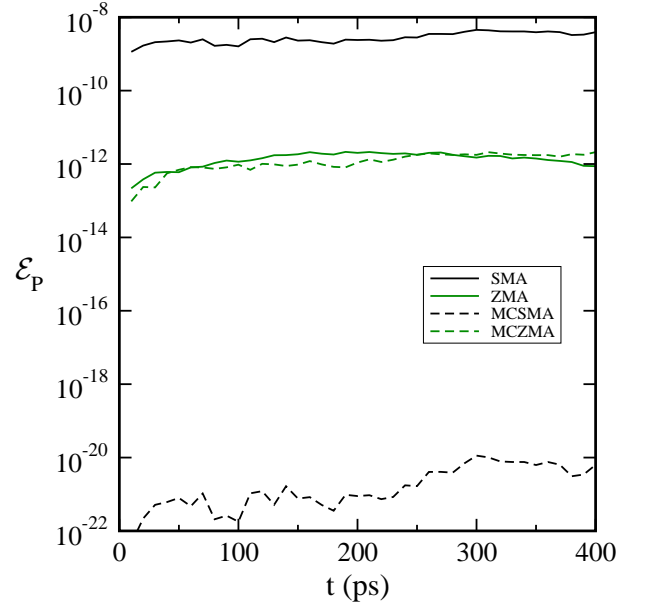}
  \caption{Solid curves show linear momentum conservation violations
    $\mathcal{E}_P$ for the large-scale approximations SMA and ZMA,
    while dashed curves show these violations for the MCSMA and MCZMA
    approximations.}
\label{f:XMA_Momentum}
\end{figure}

The solid curves in Fig.~\ref{f:XMA_Momentum} show the linear momentum
conservation violations $\mathcal{E}_P$ for the SMA and ZMA
approximations, while the dashed curves show these violations for the
MCSMA and MCZMA approximations.  We see that linear momentum violation
$\mathcal{E}_P$ are much smaller for MCSMA than the SMA approximation,
while its values are about the same for the MCZMA and ZMA
approximations.  Figure~\ref{f:XMA_Momentum} shows, however, that
linear momentum conservation is excellent for all of these large-scale
approximations.
\begin{figure}[!h]
  \includegraphics[width=.45\textwidth]{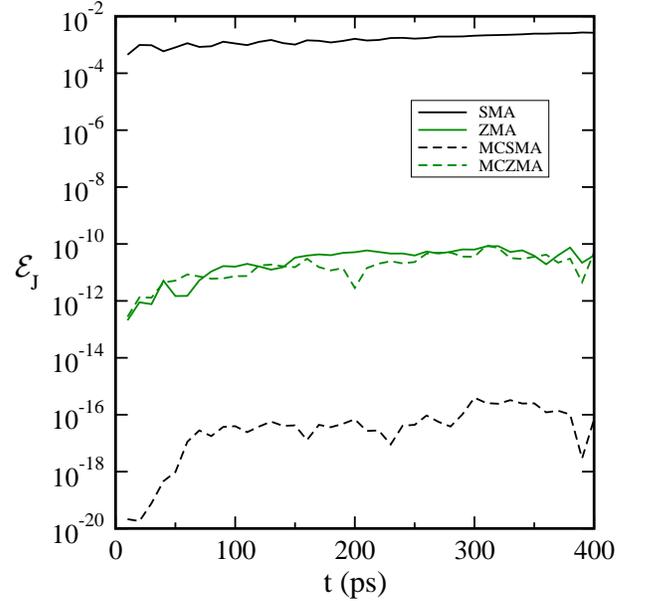}
  \caption{Solid curves show angular momentum conservation violations
    $\mathcal{E}_J$ for the large-scale approximations SMA and ZMA,
    while dashed curves show these violations for the MCSMA and MCZMA
    approximations.}
\label{f:XMA_AngularMomentum}
\end{figure}

The solid curves in Fig.~\ref{f:XMA_AngularMomentum} show the angular
momentum conservation violations $\mathcal{E}_J$ for the SMA and ZMA
approximations, while the dashed curves show these violations for the
MCSMA and MCZMA approximations.  We see that angular momentum
conservation $\mathcal{E}_J$ is much better for MCSMA than the SMA
approximation, while these violations are about the same for the MCZMA
and ZMA approximations.  Not surprisingly, the momentum conserving
approximation MCSMA has much better linear and angular momentum
conservation properties that SMA, and somewhat better momentum
conservation than the zero mode amplitude approximations MCZMA and
ZMA.

\begin{figure}[!t]
  \includegraphics[width=.45\textwidth]{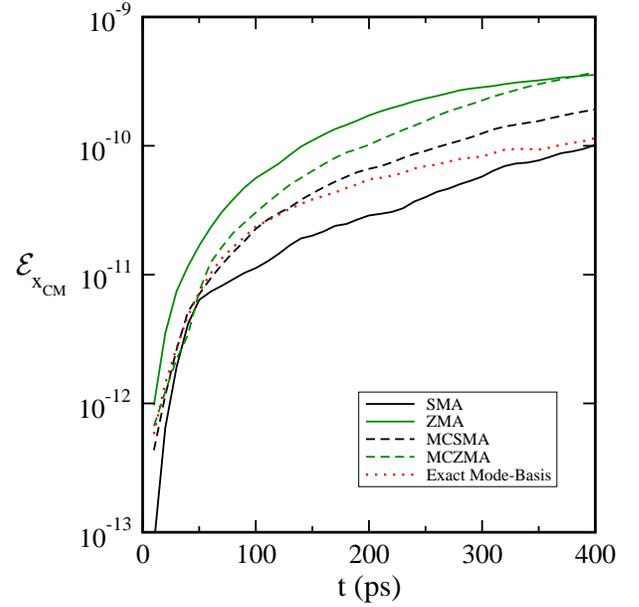}
  \caption{Solid curves show errors in the center of mass position
    $\mathcal{E}_{x_{CM}}$ for the large-scale approximations SMA and
    ZMA, while dashed curves show these errors for the MCSMA and MCZMA
    approximations.  Dotted curve shows $\mathcal{E}_{x_{CM}}$ for an
    exact mode-basis evolution with $\epsilon_\tau=10^{-12}$.}
\label{f:XMA_XCMErrors}
\end{figure}

Finally, we test how well the large-scale approximations reproduce the
macroscopic variables $\{\vec x_{CM},\vec v_{CM}, \mathbf{R},
\vec\Omega\}$ for our model problem.  We use the exact mode-basis
representation with $\epsilon_\tau=10^{-13}$ as our reference solution
in this case to test the various the large-scale approximations.  We
evaluate the differences between the approximate and the exact
solutions using the quantities $\mathcal{E}_{x_{CM}}$,
$\mathcal{E}_{v_{CM}}$, $\mathcal{E}_\Omega$ and $\mathcal{E}_q$,
defined by
\begin{eqnarray}
  \mathcal{E}_{x_{CM}}&=& \bigl|\vec x_{CM}-\vec x_{CM}^{\,Ref}\bigr|,\\
  \mathcal{E}_{v_{CM}}&=& \bigl|\vec v_{CM}-\vec v_{CM}^{\,Ref}\bigr|,\\
  \mathcal{E}_\Omega&=& \frac{\bigl|\vec\Omega-\vec\Omega^{\,Ref}\bigr|}
          {\bigl| \vec \Omega^{\,Ref}\bigr|},\\
  \mathcal{E}_q &=& \half\sqrt{\sum_i \bigl(q_i-q_i^{\,Ref}\bigr)^2}.
\end{eqnarray}
The solid curves in Fig.~\ref{f:XMA_XCMErrors} show errors in the
center of mass position, $\mathcal{E}_{x_{CM}}$, for the SMA and ZMA
approximations, and the dashed curves show these errors for the MCSMA
and MCZMA approximations.  For comparison the dotted curve in
Fig.~\ref{f:XMA_XCMErrors} shows the error in the somewhat lower
resolution, $\epsilon_\tau=10^{-12}$, evolution of the exact
mode-basis representation compared to the reference solution.  All of
these errors are very small, and are only growing slowly with time,
approximately like $\mathcal{E}_{x_{CM}}\propto t^2$. These graphs
confirm that all the large-scale approximations are able to determine
$\vec x_{CM}$ with excellent precision.
\begin{figure}[!t]
  \includegraphics[width=.45\textwidth]{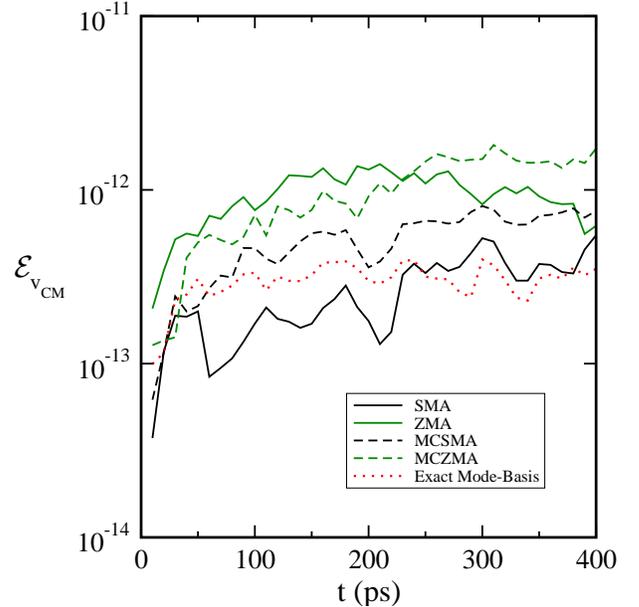}
  \caption{Solid curves show errors in the center of mass velocity
    $\mathcal{E}_{v_{CM}}$ for the large-scale approximations SMA and
    ZMA, while dashed curves show these errors for the MCSMA and MCZMA
    approximations.  Dotted curve shows $\mathcal{E}_{v_{CM}}$ for an
    exact mode-basis evolution with $\epsilon_\tau=10^{-12}$.}
\label{f:XMA_VCMErrors}
\end{figure}
The solid curves in Fig.~\ref{f:XMA_VCMErrors} show the errors in the
velocity of the center of mass, $\mathcal{E}_{v_{CM}}$, for the SMA
and ZMA approximations, and the dashed curves show these errors for
the MCSMA and MCZMA approximations.  For comparison the dotted curve
in Fig.~\ref{f:XMA_VCMErrors} shows the error in the somewhat lower
resolution, $\epsilon_\tau=10^{-12}$, evolution of the exact
mode-basis representation compared to the reference solution.  All of
these errors are very small, and appear almost constant in time at
late times. These graphs confirm that all the large-scale
approximations are able to determine $\vec v_{CM}$ with excellent
precision.

\begin{figure}[!t]
  \includegraphics[width=.45\textwidth]{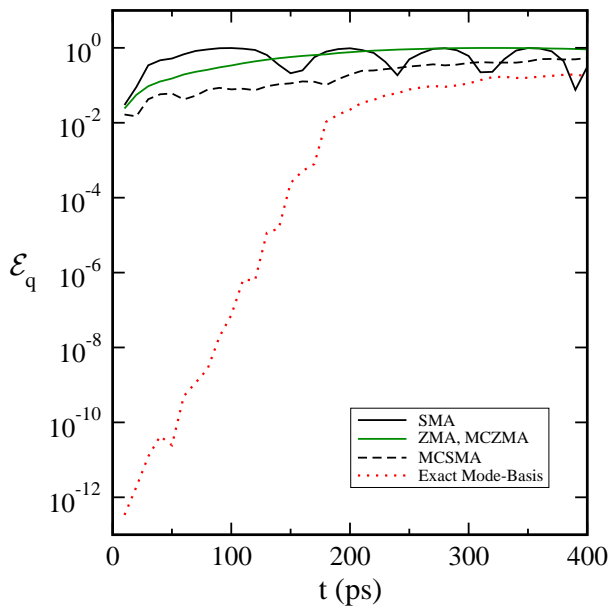}
  \caption{Solid curves show errors in the orientation parameters
    $\mathcal{E}_q$ for the large-scale approximations SMA and ZMA,
    while dashed curve shows these errors for the MCSMA approximation.
    The MCZMA errors are indistinguishable from the ZMA errors in this
    graph.  Dotted curve shows $\mathcal{E}_q$ for an exact
    mode-basis evolution with $\epsilon_\tau=10^{-12}$.  }
\label{f:XMA_qErrors}
\end{figure}

The solid curves in Fig.~\ref{f:XMA_qErrors} show the errors in the
orientation matrix $\mathbf{R}$, as measured by $\mathcal{E}_{q}$, for
the SMA and ZMA approximations, and the dashed curve shows these
errors for the MCSMA approximation.  The $\mathcal{E}_{q}$ curve for
the MCZMA approximation is indistinguishable from the ZMA curve.  For
comparison the dotted curve in Fig.~\ref{f:XMA_qErrors} shows the
error in the somewhat lower resolution, $\epsilon_\tau=10^{-12}$,
evolution of the exact mode-basis representation compared to the
reference solution.  The exact dotted curve in
Fig.~\ref{f:XMA_qErrors} shows the exponential growth at early times
which signals the presence of chaotic dynamics in this variable.  The
large scale approximations all have errors $\mathcal{E}_{q}$ that are
comparable to the late time behavior of the exact MD simulation.
While it may be surprising that the orientation of the molecule cannot
be predicted accurately from initial conditions of the molecule
using an exact MD simulation, it is not surprising in this case that
the large-scale approximations all give rather poor results for
$\mathbf{R}$ as well.

\begin{figure}[!ht]
  \includegraphics[width=.45\textwidth]{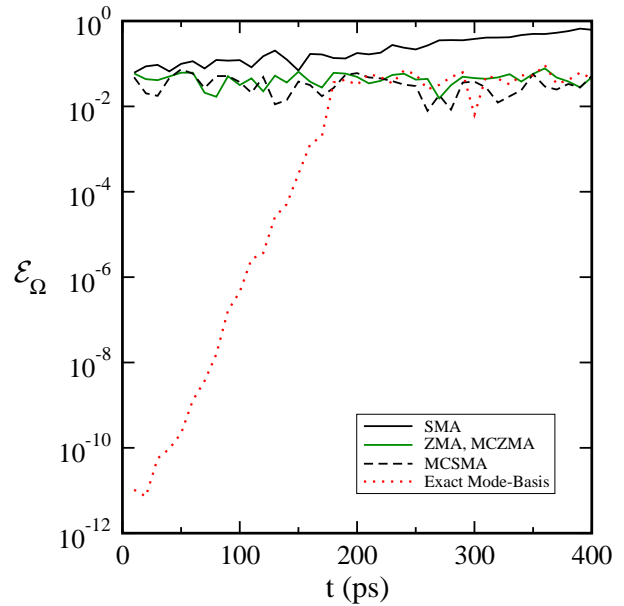}
  \caption{Solid curves show errors in the angular velocity
    $\mathcal{E}_\Omega$ for the large-scale approximations SMA and
    ZMA, while dashed curve shows these errors for the MCSMA
    approximation.  The MCZMA errors are indistinguishable from the
    ZMA errors in this graph. Dotted curve shows
    $\mathcal{E}_\Omega$ for an exact mode-basis evolution with
    $\epsilon_\tau=10^{-12}$.}
\label{f:XMA_OmegaErrors}
\end{figure}

The solid curves in Fig.~\ref{f:XMA_OmegaErrors} show errors in $\vec
\Omega$, as measured by $\mathcal{E}_{\Omega}$, for the SMA and ZMA
approximations, and the dashed curves show these errors for the MCSMA
approximation.  The MCZMA curve for $\mathcal{E}_{\Omega}$ is
indistinguishable from the ZMA curve.  For comparison the dotted curve
in Fig.~\ref{f:XMA_OmegaErrors} shows the error in the somewhat lower
resolution, $\epsilon_\tau=10^{-12}$, evolution of the exact
mode-basis representation compared to the reference solution.  The
exact dotted curve in Fig.~\ref{f:XMA_OmegaErrors} shows the
exponential growth signaling the presence of chaotic dynamics in this
variable.  Given the chaos seen in the evolution of the orientation
matrix $\mathbf{R}$ seen in Fig.~\ref{f:XMA_qErrors}, it is not at all
surprising that similar chaotic behavior is seen in the evolution of
$\vec \Omega$.  Thus it is not surprising that the large scale
approximations all have errors $\mathcal{E}_{\Omega}$ that are
comparable to the late time behavior of the exact MD simulations.

In summary: All the large scale approximations do an excellent job of
conserving linear momentum.  All the large scale approximations except
SMA do an excellent job of conserving angular momentum.  The ZMA and
MCZMA approximations do an excellent job of energy conservation, while
the SMA and MCSMA approximations do not do so well.  All of the
large-scale approximations do excellent jobs of modeling $\vec
x_{CM}$ and $\vec v_{CM}$.  None of the large scale approximations
do a good job of modeling the macroscopic orientation variables
$\mathbf{R}$ and $\vec\Omega$.  Overall then, our results show that
the ZMA and the MCZMA approximations are more reliable approximations
of the exact MD equations than the SMA and the MCSMA approximations.

\subsection{Computational Efficiency}
\label{s:ComputationalEfficiency}

\begin{figure}[!t]
  \includegraphics[width=.45\textwidth]{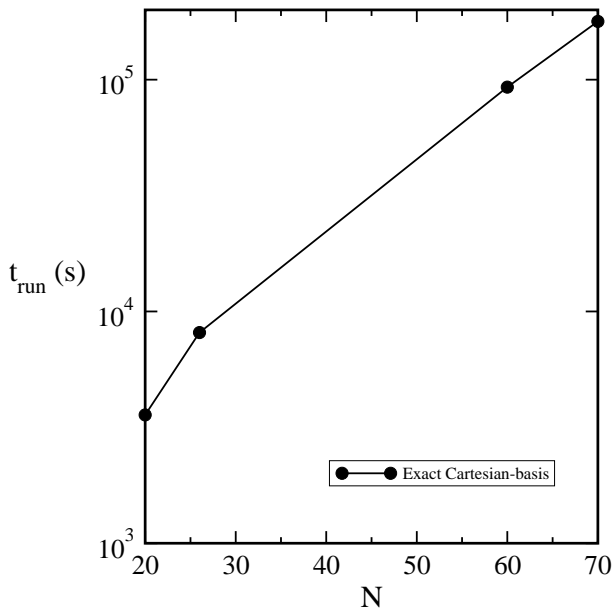}
  \caption{Total run times in seconds for the 400 ps simulations of
    the fullerene molecules $C_N$ using the Cartesian-basis version of
    the MD evolution code.  Total run time includes the runs using
    five different values of the timestep error parameters
    $\epsilon_\tau=\{10^{-6},10^{-8},10^{-10},10^{-12},10^{-13}\}$.  }
\label{f:C_RunTimes}
\end{figure}

This section briefly discusses the computational costs of the various
MD evolutions used in our tests.  Figure~\ref{f:C_RunTimes} shows the
total run times (in seconds) for the exact Cartesian-basis simulation
tests performed for each of the fullerene molecules $C_{20}$,
$C_{26}$, $C_{60}$, and $C_{70}$.  For each molecule we ran five
evolutions with $\epsilon_\tau=\{10^{-6}, 10^{-8}, 10^{-10}, 10^{-12},
10^{-13}\}$, with each evolution simulating 400 ps of the molecular
motion.  Figure~\ref{f:C_RunTimes} shows that the run times for these
tests increases exponentially as the number of atoms in the simulation
increases.  A reasonably good approximation of these total run times
is given by $t_\mathrm{run}\approx 1500\times 10^{\,N/33}$.  The code
we wrote to implement these methods was not highly optimized, so we
expect that the computational efficiency could almost certainly be
improved.

\begin{figure}[!t]
  \includegraphics[width=.45\textwidth]{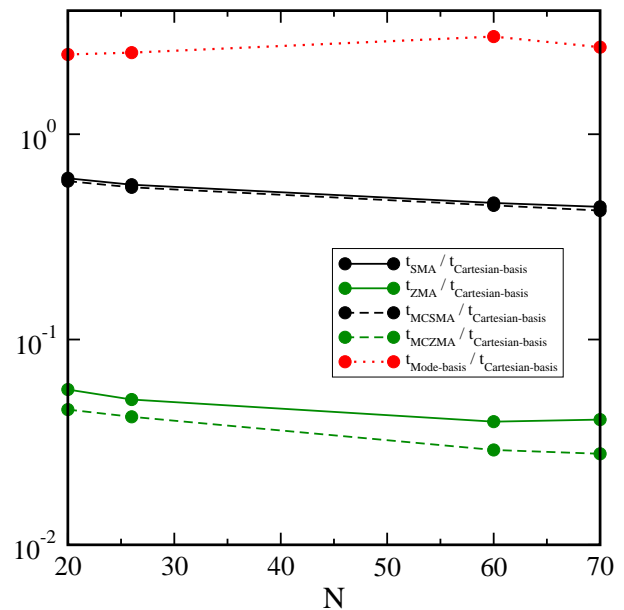}
  \caption{Relative run times for the various large-scale MD
    approximations as functions of $N$ the number of atoms in the
    simulation.  Solid curves represent the ratios of the run times of
    the SMA and ZMA approximations with the run time for the exact
    Cartesian-basis representation.  Dashed curves give the analogous
    ratios for the MCSMA and the MCZMA approximations. Dotted curve
    gives the ratio of the run times for the exact mode-basis
    representation to the exact Cartesian-basis representation.}
\label{f:Rel_Runtimes}
\end{figure}

In Fig.~\ref{f:Rel_Runtimes} we illustrate the relative computational
costs of performing evolutions using the various versions of the MD
evolution equations discussed here.  The solid curves connect the data
points that represent the ratios of the total run times for the SMA
and ZMA large-scale approximations to the total run time for the exact
Cartesian-basis evolution.  The dashed curves connect the analogous
data points for the MCSMA and the MCZMA approximations.  Finally, the
dotted curve connects the data points that represent the ratios of the
total run times for the exact mode-basis MD representation to the
total run time for exact Cartesian-basis representation.  These
results show that computations using the SMA and MCSMA approximations
are about twice as fast as those using the exact Cartesian-based
representation, while the ZMA and MCZMA approximations are more than
25 times faster.  The dotted curve in Fig.~\ref{f:Rel_Runtimes} shows
that computations using the exact mode-basis representation of MD is
two or three times slower than the exact Cartesian-basis
representation.  On the basis of computational efficiency, the MCZMA
large-scale approximation is by far the best of the various MD
simulation methods tested here.


\section{Discussion}
\label{s:Discussion}

We have developed a new mode-based representation of the classical MD
equations of motion that separate the macroscopic position and
orientation degrees of freedom of a molecule from the internal
vibrational degrees of freedom.  We have confirmed through our
numerical tests that most details of a molecular dynamical state
evolve chaotically, including the large scale orientation and angular
velocity of the molecule.  Consequently those features cannot be
predicted accurately even with exact MD simulations.  We have derived
a number of new large-scale approximations (based our new mode-based
representation) specifically designed to simulate accurately those
features of MD evolutions that are not chaotic.  We have shown through
a series of careful numerical tests that some of these large-scale
approximations give reliable, accurate predictions for the macroscopic
properties of molecular motions, including the energies, linear and
angular momentum, and the positions and velocities of their centers of
mass.  We find that one of these new approximations (MCZMA, the best
of these new large-scale approximations studied here) is more than 25
times faster than our exact Cartesian-basis MD code, while giving
comparable accuracy for the large-scale molecular properties.  We also
note that the MCZMA approximation does not depend on the mode basis
vectors $\vec e^{,\mu}_A$ at all, which makes it very easy to
implement numerically.  Thus we conclude there are many reasons to use
reliable well-tested approximations for MD simulations rather than
performing simulations using the full exact MD equations.

The tests done for this study focused on simulations of the dynamics
of single molecules.  The formalism created here has been designed,
however, to accommodate simulations of collections of molecules in a
straightforward way.  The only change that needs to be made, as we
noted earlier, is to change the single index $\scriptstyle A$ used to
identify individual atoms to a pair of indices $\scriptstyle{mA}$,
where the first index $\scriptstyle m$ determines which molecule the
particular atom belongs.  Macroscopic properties of the molecule like
the energy, $E$, linear and angular momenta, $\vec P$ and $\vec J$,
the position and velocity of the center of mass, $\vec x_{CM}$ and
$\vec v_{CM}$ should also acquire indices to identify which molecule
they belong: $\{E_m,\vec P_m,\vec J_m, \vec x_{CMm}, \vec v_{CMm}\}$.
Then, by including van der Waals and/or Coulomb
  forces in the interaction potential $U(\vec x_{mA})$, it would be
  possible to study interactions between molecules using any of the
  large-scale approximations introduced here.  The interactions
  modeled in this way should be essentially identical to those
  interactions in an exact MD simulation.

\begin{acknowledgments}
We thank Andrew McCammon and Lane Votapka for numerous helpful
discussions on molecular dynamics.  Funding and support for this work
came from the National Biomedical Computation Resource through NIH
grant P41 GM103426.  The fullerene topology graphs used in
Fig.~\ref{f:FullerineTopology} were adapted from figures produced by
R.A. Nonenmacher and published in Wikipedia under the Creative Commons
Attribution-Share Alike 3.0 Unported license.

\end{acknowledgments}

\appendix

\section{Normal Mode Basis}
\label{s:NormalModeBasis}

One natural choice for the mode-basis vectors $\vec
e^{\,\mu}_A$ are the eigenvectors of the Hessian of the potential
energy function:
\begin{equation}
0=-m_A\omega_\mu^2\,\vec e^{\,\mu}_A + \sum_B 
\frac{\partial^2 U}{\partial {\vec x_B}\partial {\vec x_A}}\cdot \vec e^{\,\mu}_B,
\label{e:ModeEquationAppendix}
\end{equation}
where $\partial^2 U/\partial {\vec x_B}\partial {\vec x_A}$ in this
equation is to be evaluated at the equilibrium state of the molecule
where $\vec x_A = \vec x_{0A}$.  Since the Hessian is a symmetric
$3N\times 3N$ dimensional real matrix, the eigenvalues $\omega_\mu^2$
and eigenvectors $\vec e^{\,\mu}_A$ are also real, and the collection
of eigenvectors form a complete basis for the $3N$ dimensional space
of vectors.  Equation~(\ref{e:ModeEquationAppendix}) is equivalent to
Newton's equation of motion, Eq.~(\ref{e:NewtonsLaw}), for the case of
very small amplitude oscillations about its equilibrium state, so we
call these $\vec e^{\,\mu}_A$ the normal-mode basis.  For stable
molecules the normal-mode frequencies are real, so all the
eigenvalues of such systems are non-negative: $\omega_\mu^2\geq 0$.

The zero-frequency modes of a molecule are represented by the
eigenvectors of the Hessian matrix of the equilibrium potential energy
function having zero eigenvalues,
\begin{eqnarray}
0=\sum_B\frac{\partial^2
U(\vec x_{0C})}{\partial \vec x_A\partial \vec x_B}\cdot \vec e^{\,\mu}_B.
\end{eqnarray}
We will now show that the eigenvectors corresponding to overall rigid
translations and rotations of the molecule are zero frequency modes.

We assume that the equilibrium state of the molecule is invariant
under rigid spatial translations and rotations.  This assumption makes
sense to ensure that the equilibrium state of interest to us is one
where the molecule is isolated and does not interact with its large
scale environment.  If the molecule is invariant under translations,
then the forces acting on the individual atoms must also be invariant.
In its equilibrium state, the total force acting on each atom must
vanish, therefore the gradient of the potential energy function must
vanish:
\begin{eqnarray}
0=\frac{\partial U(\vec x_{0C})}{\partial\vec x_A}
=\frac{\partial U(\vec x_{0C}+\lambda \vec \tau)}{\partial \vec x_A},
\label{e:TranslationEquilibriumCondtion}
\end{eqnarray}
where $\vec \tau$ is an arbitrary vector that describes the same
translation for all the atoms in the molecule, and where $\lambda$ is
an arbitrary parameter that determines the magnitude of the
translation.  We can express the forces acting on each atom of a
molecule that has been translated by an infinitesimal amount using a
Taylor series expansion:
\begin{eqnarray}
  0&=&\frac{\partial U(\vec x_{0C} + \lambda \vec \tau)}{\partial \vec x_A}
  \nonumber\\
&=& \frac{\partial U(\vec x_{0C})}{\partial \vec x_A} 
+ \lambda \sum_B \frac{\partial^2 U(\vec x_{0C})}
{\partial\vec x_{A}\partial\vec x_{B}}\cdot\vec \tau + \mathcal{O}(\lambda^2).
\label{e:TranslationExpansion}
\end{eqnarray}
The first term on the right side of Eq.~(\ref{e:TranslationExpansion})
vanishes because of the equilibrium condition,
Eq.~(\ref{e:TranslationEquilibriumCondtion}).  Therefore, the second
term on the right side of Eq.~(\ref{e:TranslationExpansion}) must also
vanish for all values of $\lambda$.  It follows that any vector
$\vec\tau$ that is the same for all the atoms in a molecule is a
zero-frequency eigenvector:
\begin{eqnarray}
0=\sum_B \frac{\partial^2 U(\vec x_{0C})}
{\partial\vec x_{A}\partial\vec x_{B}}\cdot\vec \tau.
\end{eqnarray}

The argument for the rotational invariance of the equilibrium state of
the molecule is similar.  Let
$\mathbf{R}(\lambda)$ denote a one parameter family of rotation
matrices.  We assume that $\lambda=0$ corresponds to the identity
rotation: $\mathbf{R}(0)=\mathbf{I}$. The rotational invariance of the
equilibrium state of the molecule implies that
\begin{eqnarray}
0=\frac{\partial U[\vec x_{0C}]}{\partial\vec x_A}
=\frac{\partial U[\mathbf{R}(\lambda)\cdot\vec x_{0C}]}{\partial \vec x_A}.
\label{e:RotationEquilibriumCondtion}
\end{eqnarray}
As before, we perform a Taylor expansion of the expression for the
forces acting on an equilibrium molecule that has been rotated an
infinitesimal amount:
\begin{eqnarray}0
&=&\frac{\partial U[\mathbf{R}(\lambda)\cdot\vec x_{0C}]}{\partial \vec x_A}\\
&=&\frac{\partial U[\vec x_{0C}]}{\partial \vec x_A}
+ \lambda \sum_B \frac{\partial^2 U[\vec x_{0C}]}{\partial\vec x_A
\partial \vec x_B}\cdot \left.\frac{d\mathbf{R}}{d\lambda}\right|_{\lambda=0}
\cdot \vec x_{0B}\nonumber\\
&&\qquad\qquad\qquad\qquad\qquad\qquad\qquad\quad
+\mathcal{O}(\lambda^2).
\label{e:RotationExpansion}
\end{eqnarray}
The derivative of any rotation matrix is an antisymmetric matrix.  In
this case this matrix can be written as
\begin{eqnarray}
\left.\frac{d R_{ij}}{d\lambda}\right|_{\lambda=0}
= -\sum_k\epsilon_{ijk}\theta_k,
\end{eqnarray}
where the vector $\vec\theta$ determines the direction and magnitude of
the infinitesimal rotation.  It follows that Eq.~(\ref{e:RotationExpansion})
can be re-written as
\begin{eqnarray}
0&=&\frac{\partial U[\vec x_{0C}]}{\partial \vec x_A}
+ \lambda \sum_B \frac{\partial^2 U[\vec x_{0C}]}{\partial\vec x_A
  \partial \vec x_B}\cdot \left(\vec \theta\times\vec x_{0B}\right)
+\mathcal{O}(\lambda^2).\nonumber\\
\label{e:RotationExpansionII}
\end{eqnarray}
The first term on the right side of Eq.~(\ref{e:RotationExpansionII})
vanishes because of the equilibrium condition,
Eq.~(\ref{e:RotationEquilibriumCondtion}).  Therefore, the second term
on the right side of Eq.~(\ref{e:RotationExpansionII}) must vanish for
all values of $\lambda$.  It follows that any vector of the form $\vec
\theta \times \vec x_{0A}$, where $\vec \theta$ is the same for all the atoms in 
a molecule, is a zero-frequency eigenvector:
\begin{eqnarray}
0=\sum_B \frac{\partial^2 U[\vec x_{0C}]}{\partial\vec x_A
\partial \vec x_B}\cdot \left(\vec \theta\times\vec x_{0B}\right).
\end{eqnarray}

The eigenvectors of the zero frequency modes (for generic molecules)
therefore consist of rigid translations
\begin{eqnarray}
\vec e^{\,\,t}_A(\vec\tau) = \vec \tau,
\end{eqnarray}
where $\vec \tau$ is a constant vector that determines the magnitude
and direction of the translation, and rigid rotations,
\begin{eqnarray}
\vec e^{\,\,r}_A(\vec \theta) = \vec \theta \times \vec x_{0A},
\end{eqnarray}
where $\vec \theta$ is a constant vector that determines the axis and
magnitude of the rotation.

It is easy to show from Eq.~(\ref{e:ModeEquationAppendix}) that two
eigenvectors, $\vec e^{\,\mu}_A$ and $\vec e^{\,\nu}_A$, having
different eigenvalues, $\omega_\mu^2\neq\omega_\nu^2$, are orthogonal
in the sense that,
\begin{eqnarray}
  0=\sum_A\frac{m_A}{M}\, \vec e^{\,\mu}_A \cdot \vec e^{\,\nu}_A.
\end{eqnarray}
It follows that the translations $\vec e^{\,\,t}_A(\vec \tau)=\vec
\tau$ and rotations $\vec e^{\,\,r}_A(\vec \theta)=\vec \theta\times
\vec x_{0A}$ will be orthogonal to all of the non-zero frequency mode
eigenvectors $\vec e^{\,\mu}_A$.  These orthogonality conditions are
given by:
\begin{eqnarray}
0&=& \sum_A \frac{m_A}{M} \vec e^{\,\mu}_A \cdot \vec e^{\,\,t}_A(\vec
\tau) = \sum_A \frac{m_A}{M} \vec e^{\,\mu}_A \cdot \vec \tau,\\ 0&=&
\sum_A \frac{m_A}{M} \vec e^{\,\mu}_A \cdot \vec e^{\,\,r}_A(\vec
\theta) \nonumber\\ &=& -\sum_A \frac{m_A}{M} (\vec e^{\,\mu}_A \times
\vec x_{0A})\cdot \vec \theta.
\end{eqnarray}
These orthogonality conditions will hold for arbitrary values of the
vectors $\vec \tau$ and $\vec \theta$.  Therefore these conditions can
also be written in the form
\begin{equation}
0=\sum_A \frac{m_A}{M} \vec e^{\,\mu}_A =
\sum_A \frac{m_A}{M} \vec e^{\,\mu}_A \times \vec x_{0A}.
\end{equation}
These conditions must hold for each non-zero frequency mode $\mu$, and
therefore demonstrate that the constraints on the mode-basis
eigenvectors given in Eqs.~(\ref{e:TranslationConstraint}) and
(\ref{e:RotationConstraint}) are satisfied by the normal-mode basis
vectors.

\section{Quaternion Representation of ${\mathbf R}(t)$}
\label{s:QuaternionRepresentation}

The differential equation that determines the rotation matrix $\mathbf{R}(t)$,
\begin{eqnarray}
\frac{dR_{ij}}{dt}=-\sum_{k\ell}\epsilon_{ik\ell}\,\Omega^\ell R_{kj},
\label{e:dRdtAppendix}
\end{eqnarray}
can be integrated numerically directly.  Unfortunately the
accumulation of truncation and roundoff errors in this direct approach
inevitably produces a solution that is no longer a rotation matrix,
and there is no reliable way to project out these errors to retrieve
the correct $\mathbf{R}(t)$.  A better approach is to adopt some
parametric representation of the three-dimensional space of rotation
matrices, then to convert Eq.~(\ref{e:dRdtAppendix}) into a system of
equations for the evolution of those parameters, and finally to
integrate that parametric representation numerically.  For example,
one common representation of the rotation matrices uses the Euler
angles as parameters.  Since the Euler angle representation is not one
to one (at a few singular points), a better representation uses unit
quaternions which do provide a one to one representation.  We use the
quaternion representation for our numerical solutions of
$\mathbf{R}(t)$.  Let $q_0$ represent the real part, and $q_1$, $q_2$,
and $q_3$ the three independent imaginary parts of a quaternion with
\begin{equation}
q_0^2(t)+q_1^2(t)+q_2^2(t)+q_3^2(t) = 1.
\label{e:unit}
\end{equation}
This equation defines a unit three-sphere in this parameter space, so
the space of possible parameter values is three-dimensional. A general
rotation matrix $\mathbf{R}$ can be written in terms of these
quaternion parameters in the following way,
\begin{equation}
\mathbf{R} = 2\!\begin{pmatrix} 
q_0^2+q_1^2-\tfrac{1}{2} &\,\, q_1q_2-q_0q_3 &\,\, q_1q_3+q_0q_2   \\
q_1q_2+q_0q_3  &\,\, q_0^2+q_2^2-\tfrac{1}{2} &\,\, q_2q_3-q_0q_1 \\
q_1q_3-q_0q_2   &\,\, q_2q_3+q_0q_1 &\,\,q_0^2+q_3^2-\tfrac{1}{2} \\
\end{pmatrix}\!.
\label{e:quat2mat2}
\end{equation}
It is then straightforward to transform the rotation matrix evolution
Eq.~(\ref{e:dRdtAppendix}) into an equation for the evolution of the
quaternion parameters.  The result is
\begin{align}
  \frac{dq_0(t)}{dt} &= -\tfrac{1}{2}(\Omega_x q_1 +\Omega_y q_2 + \Omega_z q_3),
\label{e:diffquat0}\\
\frac{dq_1(t)}{dt} &= \tfrac{1}{2}(\Omega_x q_0 +\Omega_y q_3 - \Omega_z q_2),
\label{e:diffquat1}\\
\frac{dq_2(t)}{dt} &= \tfrac{1}{2}(-\Omega_x q_3 +\Omega_y q_0 + \Omega_z q_1),
\label{e:diffquat2}\\
\frac{dq_3(t)}{dt} &= \tfrac{1}{2}(\Omega_x q_2 -\Omega_y q_1 + \Omega_z q_0).
\label{e:diffquat3}
\end{align}
Given a solution to these equations for $q_0(t)$, $q_1(t)$, $q_2(t)$
and $q_3(t)$, it is easy to reconstruct the rotation matrix
$\mathbf{R}(t)$ using Eq.~(\ref{e:quat2mat2}).

The constraint,
\begin{equation}
  \mathcal{C} \equiv q_0^2 + q_1^2 + q_2^2 + q_3^2 -1,
  \label{e:qconstraint}
\end{equation}
which measures how well the quaternion parameters remain on the unit
three-sphere, is preserved by the evolution defined by
Eqs.~(\ref{e:diffquat0})--(\ref{e:diffquat3}).  In
particular these evolution equations imply
\begin{equation}
  \frac{d\mathcal{C}}{dt}=0.
\end{equation}
Solving Eqs.~(\ref{e:diffquat0})--(\ref{e:diffquat3})
numerically will nevertheless result in truncation level violations of
this constraint, so it is necessary to re-scale the solution
periodically to ensure that $\mathcal{C}=0$.  Without this re-scaling
the $\mathbf{R}(t)$ constructed using Eq.~(\ref{e:quat2mat2}) will not
be a rotation matrix.

The numerical solution of
Eqs.~(\ref{e:diffquat0})--(\ref{e:diffquat3}) can be improved by
adding constraint damping terms to the system.  These extra terms
vanish whenever the constraints are satisfied, thus leaving the
desired solutions unchanged.  But these constraint damping terms are
chosen to drive the solution back toward the constraint satisfying
surface whenever small numerical constraint violations inevitably
occur.  The resulting quaternion evolution equation with constraint
damping is given by,
\begin{align}
\frac{dq_0(t)}{dt} &= -\tfrac{1}{2}(\Omega_x q_1 +\Omega_y q_2 + \Omega_z q_3)
-\tfrac{1}{8}\,\eta\, q_0\, \mathcal{C},
\label{e:diffquat0c}\\
\frac{dq_1(t)}{dt} &= \tfrac{1}{2}(\Omega_x q_0 +\Omega_y q_3 - \Omega_z q_2)
-\tfrac{1}{8}\,\eta\, q_1\, \mathcal{C},
\label{e:diffquat1c}\\
\frac{dq_2(t)}{dt} &= \tfrac{1}{2}(-\Omega_x q_3 +\Omega_y q_0 + \Omega_z q_1)
-\tfrac{1}{8}\,\eta\, q_2\, \mathcal{C},
\label{e:diffquat2c}\\
\frac{dq_3(t)}{dt} &= \tfrac{1}{2}(\Omega_x q_2 -\Omega_y q_1 + \Omega_z q_0)
-\tfrac{1}{8}\,\eta\, q_3\, \mathcal{C}.
\label{e:diffquat3c}
\end{align}
These equations imply the following evolution equation for the
constraint,
\begin{equation}
\frac{d\mathcal{C}}{dt} = -\eta (\mathcal{C}+1)\mathcal{C},
\end{equation}
which drives constraint violations $\mathcal{C}$ toward zero
exponentially on a timescale set by the constant $\eta$ when $\eta>0$.
This is the form of the rotation matrix evolution equations used in
all of our numerical solutions of the various representations of the MD
equations.

\section*{References}
\label{s:References}
\bibliography{mybibfile}

\end{document}